\documentclass{amsart}

\usepackage{amsmath}
\usepackage{amsfonts}
\usepackage{amssymb}

\theoremstyle{proclaim}
\newtheorem{theorem}{Theorem}[section]
\newtheorem{corollary}[theorem]{Corollary}
\newtheorem{proposition}[theorem]{Proposition}
\newtheorem{lemma}[theorem]{Lemma}

\theoremstyle{definition}
\newtheorem{definition}[theorem]{Definition}
\newtheorem{example}[theorem]{Example}
\newtheorem*{remarks}{Remarks on Definition 3.2}
\newtheorem*{remark}{Remark}
\newtheorem*{acknowledgment}{Acknowledgment}

\begin{document}

\date{2010-6-19}

\title[The Szemer\'{e}di property]{The Szemer\'{e}di property
in ergodic W*-dynamical systems}

\author{Conrad Beyers, Rocco Duvenhage and Anton Str\"{o}h}

\address{CONRAD BEYERS, Department of Mathematics and Applied Mathematics,
University of Pretoria, Pretoria, 0002, South Africa}

\email{conrad.beyers@up.ac.za}

\address{ROCCO DUVENHAGE, Department of Mathematics and Applied
Mathematics, (Current address: Department of Physics), University of
Pretoria, Pretoria, 0002, South Africa}

\email{rocco.duvenhage@up.ac.za}

\address{ANTON STR\"{O}H, Department of Mathematics and Applied Mathematics,
University of Pretoria, Pretoria, 0002, South Africa}

\email{anton.stroh@up.ac.za}

\begin{abstract}
We study weak mixing of all orders for asymptotically abelian weakly
mixing state preserving C*-dynamical systems, where the dynamics is
given by the action of an abelian second countable locally compact
group which contains a F\o lner sequence satisfying the Tempelman
condition. For a smaller class of groups (which include
$\mathbb{Z}^{q}$\ and $\mathbb{R}^{q}$) this is then used to show
that an asymptotically abelian ergodic W*-dynamical system either
has the ``Szemer\'{e}di property'' or contains a nontrivial
subsystem (a ``compact factor'') that does. A van der Corput lemma
for Hilbert space valued functions on the group is one of our main
technical tools.
\end{abstract}



\maketitle

\section{Introduction}

Consider a measure preserving dynamical system $(X,\Sigma,\nu,T)$, i.e. $T$ is
an invertible measure preserving transformation of a probability space
$(X,\Sigma,\nu)$, namely a set $X$ with $\sigma$-algebra $\Sigma$ on which
$\nu$ is a measure with $\nu(X)=1$. Furstenberg \cite{F1}, \cite{F2} proved
that for any such system
\[
\liminf_{n\rightarrow\infty}\frac{1}{N}\sum_{n=1}^{N}\nu(A\cap T^{-n}(A)\cap
T^{-2n}(A)\cap...\cap T^{-kn}(A))>0
\]
if $\nu(A)>0$. We will refer to this as a \textit{Szemer\'{e}di property};
also see \cite{FKO}. To prove this result requires (among other things) a
structure theory in terms of so-called weakly mixing systems and compact systems.

Weak mixing is an important notion in ergodic theory, introduced by Koopman
and von Neumann \cite{KvN} in 1932 for actions of the group $\mathbb{R}$.
Iterates of $T$ above can be viewed as an action of the group $\mathbb{Z}$,
and in this case the system above is called weakly mixing if
\[
\lim_{N\rightarrow\infty}\frac{1}{N}\sum_{n=1}^{N}\left|  \nu(A\cap
T^{-n}(B))-\nu(A)\nu(B)\right|  =0
\]
for all $A,B\in\Sigma$. Under this assumption Furstenberg proved that the
system is in fact weakly mixing of all orders, namely
\begin{equation}
\lim_{N\rightarrow\infty}\frac{1}{N}\sum_{n=1}^{N}\left|  \nu(A_{0}\cap
T^{-m_{1}n}(A_{1})\cap...\cap T^{-m_{k}n}(A_{k}))-\nu(A_{0})\nu(A_{1}%
)...\nu(A_{k})\right|  =0 \tag{1.1}%
\end{equation}
for all $A_{0},...,A_{k}\in\Sigma$, all
$m_{1},...,m_{k}\in\mathbb{N}$ with $m_{1}<m_{2}<...<m_{k}$, and all
$k\in\mathbb{N=}\left\{  1,2,3,...\right\} $, from which the
Szemer\'{e}di property then follows easily for weakly mixing
systems. On the other hand, the system is called compact if the
orbit $\left\{  f\circ T^{n}:n\in\mathbb{Z}\right\}  $ of every
$f\in L^{2}(\nu)$ is relatively compact in $L^{2}(\nu)$. Such
systems can also be shown to have the Szemer\'{e}di property. As one
might expect, these and related ideas have been studied for actions
of more general groups; see for example \cite{Dye}, \cite{B2}
(Section 4), \cite{BMZ} and \cite{B3}.

The results of this paper form part of a programme to extend the structure
theorems developed by Furstenberg and others to the operator algebraic
setting. In this paper we first study weak mixing of all orders in a
non-commutative C*-algebraic setting where $(X,\Sigma,\nu,T)$ is replaced by a
C*-dynamical system $(A,\omega,\tau)$ where $\omega$ is a state on the unital
C*-algebra $A$, and $\tau$ a group of $\ast$-automorphisms of $A$ keeping
$\omega$ invariant. This problem has also been studied by Niculescu,
Str\"{o}h, and Zsid\'{o} \cite{NSZ} for actions of $\mathbb{Z}$, however we
allow more general groups, namely abelian second countable locally compact
groups which contains a F\o lner sequence satisfying certain conditions. The
role of a F\o lner sequence is to replace the sequence of sets $\{1,...,n\}$
appearing in the averages in the expressions above. We then proceed to the
Szemer\'{e}di property for compact C*-dynamical systems (Section 5). Finally,
in Section 6, we use the results of the previous sections to study ergodic
W*-dynamical systems (where $A$ above is a $\sigma$-finite von Neumann
algebra), however we only show that an asymptotic abelian ergodic system
either has the Szemer\'{e}di property, or has a subsystem (called a factor)
that has this property . This final result (Theorem 6.10) is proved for a
smaller class of groups which however still contains $\mathbb{Z}^{q}$\ and
$\mathbb{R}^{q}$. Many of the intermediate results hold for more general
groups or semigroups, as we will indicate. The asymptotic abelianness we refer
to here is of a relatively weak form, namely ``in the average'' or ``in
density'' as defined in Section 4. Asymptotic abelianness is needed to handle
the weakly mixing case, while the compact case works without it, however in
the latter we assume $\omega$ to be tracial while in the former we do not. So
a certain level of commutativity is always present.

The largest part of the paper is devoted to the weakly mixing case. One of the
technical tools we use in this case is a so-called van der Corput lemma which
we discuss in Section 2. This type of lemma and related inequalities, inspired
by the classical van der Corput difference theorem and van der Corput
inequality, have been used by Bergelson et al \cite{B1}, \cite{BMZ},
Furstenberg \cite{F3}, Niculescu, Str\"{o}h, and Zsid\'{o} \cite{NSZ}, and
others, to study polynomial ergodic theorems, nonconventional ergodic
averages, and noncommutative recurrence, for example. We extend the van der
Corput lemma to more general groups, namely second countable amenable locally
compact groups. The main result of this section is given by Theorem 2.6. After
some preliminaries on weak mixing in Section 3, we devote Section 4 to showing
how weak mixing implies weak mixing of all orders. The form of weak mixing of
all orders we prove, involves replacing the multiplication with $m_{1}%
,...,m_{k}$ in (1.1), by homomorphisms of the group over which we work, and
this motivates why we incorporate such homomorphisms in a
generalized\ definition of weak mixing in Section 3. The main result of
Section 4 is Theorem 4.6.

\section{A van der Corput lemma}

This section is devoted to proving a van der Corput lemma, stated in Theorem
2.6. Our proof of the van der Corput lemma will roughly follow that of
\cite{F3} over the group $\mathbb{Z}$.

For $(Y,\mu)$ a measure space and $\mathfrak{H}$ a Hilbert space,
consider a bounded function $f:\Lambda\rightarrow\mathfrak{H}$ with
$\Lambda\subset Y$ measurable and $\mu(\Lambda)<\infty$, and
$\left\langle f(\cdot),x\right\rangle $ measurable for every
$x\in\mathfrak{H}$. Define $\int_{\Lambda}fd\mu\in \mathfrak{H}$ by
requiring
\begin{equation}
\left\langle \int_{\Lambda}fd\mu,x\right\rangle :=\int_{\Lambda}\left\langle
f(y),x\right\rangle d\mu(y) \tag{2.1}%
\end{equation}
for all $x\in\mathfrak{H}$. We will often use the notation $\int_{\Lambda
}f(y)dy=\int_{\Lambda}fd\mu$, since there will be no ambiguity in the measure
being used. Iterated integrals (when they exist) will be written as $\int
_{B}\int_{A}f(y,z)dydz$, which of course simply means $\int_{B}\left[
\int_{A}f(y,z)dy\right]  dz$, and similarly for triple integrals.

In a group $G$ we will use the notations $Vg:=\{vg:v\in V\},~VW:=\{vw:v\in
V,w\in W\},~V^{-1}:=\{v^{-1}:v\in V\}$, etc. for any $V,W\subset G$ and $g\in
G$, and we will use multiplicative notation even when working in an abelian group.

In Sections 2 to 4 of this paper $G$ denotes an abelian second countable
locally compact group with identity $e$, and regular Haar measure $\mu$. Since
$G$ is abelian, it is amenable. In this section and the next the abelianness
of $G$ is in fact not crucial; the proofs go through even if $G$ is not
abelian but still amenable, and $\mu$ is right invariant (but see the remarks
just before Theorem 2.6). Unfortunately in Section 4 this is not the case.

Since $G$ is second countable and locally compact, it is $\sigma$-compact and
hence its amenability (even for a nonabelian group) is equivalent to the
existence of a F\o lner sequence $\left(  \Lambda_{n}\right)  $ in $G$ defined
as follows:

\begin{definition} 
A \emph{F\o lner sequence} in $G$ is a sequence $\left(
\Lambda_{n}\right)  $ of compact subsets of $G$ such that
$0<\mu(\Lambda_{n})$ for all $n$, and
\begin{equation}
\lim_{n\rightarrow\infty}\frac{\mu\left(  \Lambda_{n}\Delta(\Lambda
_{n}g)\right)  }{\mu(\Lambda_{n})}=0 \tag{2.1.1}%
\end{equation}
for all $g\in G$.
\end{definition}

Refer to Theorem 4 in \cite{E1} and Theorems 1 and 2 in \cite{E2} for a very
clear exposition of this. In fact, these papers show that we can choose a F\o
lner sequence with stronger properties than those in Definition 2.1, but our
definition will suffice for this paper. Furthermore, Theorem 3 in \cite{E1}
shows that Definition 2.1 implies uniform convergence of (2.1.1) on compact
sets, ie.
\[
\lim_{n\rightarrow\infty}\sup_{g\in K}\frac{\mu\left(  \Lambda_{n}%
\Delta(\Lambda_{n}g)\right)  }{\mu(\Lambda_{n})}=0
\]
for any non-empty compact $K\subset G$. We will have occasion to use this
important fact later on. Throughout Sections 2 to 4, $\left(  \Lambda
_{n}\right)  $ will denote a F\o lner sequence in $G$. At the end of Section
4, we briefly consider simple examples of such sequences in $\mathbb{Z}^{q}$
and $\mathbb{R}^{q}$.

We assume second countability, since for second countable topological spaces
$X,Y$, and their Borel $\sigma$-algebras $S,T$, the product $\sigma$-algebra
obtained from $S,T$ is the same as the Borel $\sigma$-algebra of the
topological space $X\times Y$. This is needed in order to apply Fubini's
theorem, which requires measurability in the product $\sigma$-algebra.

\begin{proposition} 
Consider a bounded $f:G\rightarrow \mathfrak{H}$
with $\mathfrak{H}$ a Hilbert space, such that $\left\langle
f(\cdot),x\right\rangle $ is Borel measurable for every
$x\in\mathfrak{H}$. Then
\[
\lim_{m\rightarrow\infty}\left|  \left|  \frac{1}{\mu(\Lambda_{m})}
\int_{\Lambda_{m}}fd\mu-\frac{1}{\mu(\Lambda_{m})}\frac{1}{\mu(\Lambda_{n}
)}\int_{\Lambda_{m}}\int_{\Lambda_{n}}f(gh)dhdg\right|  \right|  =0
\]
for every $n$.
\end{proposition}

\begin{proof}
By (2.1) and Fubini's theorem
\[
\int_{\Lambda_{m}}\int_{\Lambda_{n}}f(gh)dhdg=\int_{\Lambda_{n}}\int
_{\Lambda_{m}}f(gh)dgdh
\]
and in particular these iterated integrals exists. From this and the fact that
$\mu$ is a right invariant measure, we have
\begin{align*}
&  \left|  \left|  \frac{1}{\mu(\Lambda_{m})}\int_{\Lambda_{m}}fd\mu
-\frac{1}{\mu(\Lambda_{m})}\frac{1}{\mu(\Lambda_{n})}\int_{\Lambda_{m}}%
\int_{\Lambda_{n}}f(gh)dhdg\right|  \right| \\
&  =\left|  \left|  \frac{1}{\mu(\Lambda_{n})}\frac{1}{\mu(\Lambda_{m})}%
\int_{\Lambda_{n}}\left[  \int_{\Lambda_{m}}f(g)dg-\int_{\Lambda_{m}%
}f(gh)dg\right]  dh\right|  \right| \\
&  =\left|  \left|  \frac{1}{\mu(\Lambda_{n})}\frac{1}{\mu(\Lambda_{m})}%
\int_{\Lambda_{n}}\left[  \int_{\Lambda_{m}}f(g)dg-\int_{\Lambda_{m}%
h}f(g)dg\right]  dh\right|  \right| \\
&  =\left|  \left|  \frac{1}{\mu(\Lambda_{n})}\frac{1}{\mu(\Lambda_{m})}%
\int_{\Lambda_{n}}\left[  \int_{\Lambda_{m}\backslash\left(  \Lambda_{m}%
\cap(\Lambda_{m}h)\right)  }f(g)dg-\int_{\left(  \Lambda_{m}h\right)
\backslash\left(  \Lambda_{m}\cap(\Lambda_{m}h)\right)  }f(g)dg\right]
dh\right|  \right|
\end{align*}
But if$\ b\in\mathbb{R}$ is an upper bound for $\left|  \left|
f(G)\right| \right|  $ , we have
\[
\left|  \left|  \int_{\Lambda_{m}\backslash\left(  \Lambda_{m}\cap(\Lambda
_{m}h)\right)  }f(g)dg-\int_{\left(  \Lambda_{m}h\right)  \backslash\left(
\Lambda_{m}\cap(\Lambda_{m}h)\right)  }f(g)dg\right|  \right|  \leq
b\sup_{h\in\Lambda_{n}}\mu\left(  \Lambda_{m}\Delta(\Lambda_{m}h)\right)
\]
therefore
\begin{align*}
&  \left|  \left|  \frac{1}{\mu(\Lambda_{m})}\int_{\Lambda_{m}}fd\mu
-\frac{1}{\mu(\Lambda_{m})}\frac{1}{\mu(\Lambda_{n})}\int_{\Lambda_{m}}%
\int_{\Lambda_{n}}f(gh)dhdg\right|  \right| \\
&  \leq\frac{1}{\mu(\Lambda_{m})}b\sup_{h\in\Lambda_{n}}\mu\left(  \Lambda
_{m}\Delta(\Lambda_{m}h)\right) \\
&  \rightarrow0
\end{align*}
as $m\rightarrow\infty$.
\end{proof}

\begin{lemma} 
Let $\mathfrak{H}$ be a Hilbert space, $(Y,\mu)$ a measure space,
and $\Lambda\subset Y$ a measurable set with $\mu(\Lambda)<\infty$.
Consider an $f:\Lambda\rightarrow\mathfrak{H}$ with $\left|  \left|
f(\cdot)\right|  \right|  $ measurable, and $\left\langle
f(\cdot),x\right\rangle $ measurable for every $x\in \mathfrak{H}$,
and with $\int_{\Lambda}\left|  \left|  f(y)\right|  \right|
dy<\infty$ (which means $\int_{\Lambda}fd\mu$ exists). Then
\[
\left\|  \int_{\Lambda}fd\mu\right\|  ^{2}\leq\mu(\Lambda)\int_{\Lambda
}\left|  \left|  f(y)\right|  \right|  ^{2}dy
\]
\end{lemma}

\begin{proof}
From the definition of $\int_{\Lambda}fd\mu$,
\begin{align*}
\left\|  \int_{\Lambda}fd\mu\right\|  ^{2}  &  =\int_{\Lambda}\left[
\int_{\Lambda}\operatorname{Re}\left\langle f(y),f(z)\right\rangle dz\right]
dy\\
&  \leq\frac{1}{2}\int_{\Lambda}\left[  \int_{\Lambda}\left(  \left|
\left| f(y)\right|  \right|  ^{2}+\left|  \left|  f(z)\right|
\right|  ^{2}\right) dz\right]  dy\text{.}
\end{align*}
\end{proof}

\begin{proposition} 
Consider the situation in Proposition 2.2. Assume furthermore that
$F:G\times G\rightarrow \mathbb{C}:(g,h)\mapsto\left\langle
f(g),f(h)\right\rangle $ is Borel measurable, and that
$\Lambda_{1},\Lambda_{2}\subset G$ are Borel sets with
$\mu(\Lambda_{j})<\infty$. Then
\[
\left|  \left|  \int_{\Lambda_{2}}\int_{\Lambda_{1}}f(gh)dhdg\right|
\right|
^{2}\leq\mu(\Lambda_{2})\int_{\Lambda_{1}}\int_{\Lambda_{1}}\int_{\Lambda_{2}
}\left\langle f(gh_{1}),f(gh_{2})\right\rangle dgdh_{1}dh_{2}
\]
and in particular these integrals exist.
\end{proposition}

\begin{proof}
The double integral exists as in Proposition 2.2's proof. Let's now
consider the triple integral. Since $F$ is Borel measurable and $G$
's product is continuous, $(g,h_{1})\mapsto\left\langle f(gh_{1}
),f(gh_{2})\right\rangle $ is Borel measurable on $G\times G=G^{2}$
and hence measurable in the product $\sigma$-algebra on $G^{2}$. By
Fubini's theorem we have
\[
\int_{\Lambda_{1}}\int_{\Lambda_{2}}\left\langle f(gh_{1}),f(gh_{2}%
)\right\rangle dgdh_{1}=\int_{\Lambda_{1}\times\Lambda_{2}}\left\langle
f(gh_{1}),f(gh_{2})\right\rangle d(h_{1},g)
\]
and in particular the iterated integral exists. Furthermore, $G\times
G^{2}\rightarrow G^{2}:(h_{2},h_{1},g)\mapsto(gh_{1},gh_{2})$ is continuous,
so $G\times G^{2}\rightarrow\mathbb{C}:(h_{2},h_{1},g)\mapsto\left\langle
f(gh_{1}),f(gh_{2})\right\rangle $ is measurable in the product $\sigma
$-algebra of $G$ and $G^{2}$. Hence by Fubini's theorem
\[
\int_{\Lambda_{1}}\int_{\Lambda_{1}\times\Lambda_{2}}\left\langle
f(gh_{1}),f(gh_{2})\right\rangle d(h_{1},g)dh_{2}=\int_{\Lambda_{1}%
\times\Lambda_{1}\times\Lambda_{2}}\left\langle f(gh_{1}),f(gh_{2}%
)\right\rangle d(h_{2},h_{1},g)
\]
and in particular, the triple integral exists, and we can do the three
integrals in any order. By Lemma 2.3 it follows that
\begin{align*}
&  \left|  \left|  \int_{\Lambda_{2}}\int_{\Lambda_{1}}f(gh)dhdg\right|
\right|  ^{2}\\
&  \leq\mu(\Lambda_{2})\int_{\Lambda_{2}}\left|  \left|  \int_{\Lambda_{1}%
}f(gh)dh\right|  \right|  ^{2}dg\\
&  =\mu(\Lambda_{2})\int_{\Lambda_{1}}\int_{\Lambda_{1}}\int_{\Lambda_{2}%
}\left\langle f(gh_{1}),f(gh_{2})\right\rangle dgdh_{1}dh_{2}%
\end{align*}
and note in particular that the last equality proves that $g\mapsto\left|
\left|  \int_{\Lambda_{1}}f(gh)dh\right|  \right|  ^{2}$ is measurable (and
therefore its square root too), which means that Lemma 2.3 does indeed apply
to this situation.
\end{proof}

\begin{proposition} 
Consider the situation in Proposition 2.2. Assume that $F:G\times
G\rightarrow\mathbb{C}:(g,h)\mapsto \left\langle
f(g),f(h)\right\rangle $ is Borel measurable. Then $\int
_{\Lambda}\left\langle f(g),f(gh)\right\rangle dg$ exists for all
measurable $\Lambda\subset G$ with $\mu(\Lambda)<\infty$, and all
$h\in G$. Assume that
\[
\gamma_{h}:=\lim_{n\rightarrow\infty}\frac{1}{\mu(\Lambda_{n})}\int
_{\Lambda_{n}}\left\langle f(g),f(gh)\right\rangle dg
\]
exists for all $h\in G$. Then
\[
\lim_{m\rightarrow\infty}\frac{1}{\mu(\Lambda_{m})}\int_{\Lambda_{n}}
\int_{\Lambda_{n}}\int_{\Lambda_{m}}\left\langle f(gh_{1}),f(gh_{2}
)\right\rangle dgdh_{1}dh_{2}=\int_{\Lambda_{n}}\int_{\Lambda_{n}}
\gamma_{h_{1}^{-1}h_{2}}dh_{1}dh_{2}
\]
for all $n$, and in particular these integrals exist.
\end{proposition}

\begin{proof}
The triple integral exists by Proposition 2.4. Let $b$
be an upper bound for $(g,h)\mapsto\left|  \left\langle
f(g),f(h)\right\rangle \right|  $, which exists since $f$ is
bounded. Fix any $m\in\mathbb{N}$, and set
\[
A_{n}(h_{1},h_{2}):=\frac{1}{\mu(\Lambda_{n})}\int_{\Lambda_{n}}\left\langle
f(gh_{1}),f(gh_{2})\right\rangle dg
\]
for all $h_{1},h_{2}\in\Lambda_{m}$ and all $n$. Note that
$A_{n}(h_{1} ,h_{2})$ exists and is a measurable function of $\left(
h_{1},h_{2}\right)  $ by Fubini's Theorem. Since $F$ is Borel, and
$G\rightarrow G^{2} :g\mapsto(g,gh)$ is continuous, the map
$G\rightarrow\mathbb{C}:g\mapsto \left\langle
f(g),f(gh)\right\rangle $ is Borel for every $h\in G$. Now,
\begin{align*}
\left|  A_{n}(h_{1},h_{2})-\gamma_{h_{1}^{-1}h_{2}}\right|   &  \leq
\frac{1}{\mu(\Lambda_{n})}\left|  \int_{\Lambda_{n}}\left\langle
f(gh_{1}),f(gh_{2})\right\rangle dg-\int_{\Lambda_{n}}\left\langle
f(g),f(gh_{1}^{-1}h_{2})\right\rangle dg\right| \\
&  +\left|  \frac{1}{\mu(\Lambda_{n})}\int_{\Lambda_{n}}\left\langle
f(g),f(gh_{1}^{-1}h_{2})\right\rangle dg-\gamma_{h_{1}^{-1}h_{2}}\right|
\end{align*}
for all $h_{1}\in G$ and $h_{2}\in G$. But since $\mu$ is a right
invariant measure
\begin{align*}
&  \frac{1}{\mu(\Lambda_{n})}\left|  \int_{\Lambda_{n}}\left\langle
f(gh_{1}),f(gh_{2})\right\rangle dg-\int_{\Lambda_{n}}\left\langle
f(g),f(gh_{1}^{-1}h_{2})\right\rangle dg\right| \\
&  =\frac{1}{\mu(\Lambda_{n})}\left|  \int_{\Lambda_{n}h_{1}}\left\langle
f(g),f(gh_{1}^{-1}h_{2})\right\rangle dg-\int_{\Lambda_{n}}\left\langle
f(g),f(gh_{1}^{-1}h_{2})\right\rangle dg\right| \\
&  \leq\frac{1}{\mu(\Lambda_{n})}\left[  \int_{(\Lambda_{n}h_{1}
)\backslash\Lambda_{n}}\left|  \left\langle f(g),f(gh_{1}^{-1}h_{2}
)\right\rangle \right|
dg+\int_{\Lambda_{n}\backslash(\Lambda_{n}h_{1} )}\left|
\left\langle f(g),f(gh_{1}^{-1}h_{2})\right\rangle \right|
dg\right] \\
&  \leq\frac{\mu\left(  \Lambda_{n}\Delta(\Lambda_{n}h_{1})\right)
} {\mu(\Lambda_{n})}b
\end{align*}
for all $h_{1}\in G$. Hence $\lim_{m\rightarrow\infty}A_{n}(h_{1}
,h_{2})=\gamma_{h_{1}^{-1}h_{2}}$.

Furthermore, $\left|  A_{n}(h_{1},h_{2})\right|  \leq\frac{1}{\mu(\Lambda
_{n})}\int_{\Lambda_{n}}bdg=b$, which implies that the sequence $A_{n}$ is
dominated by $B:\Lambda_{m}\times\Lambda_{m}\rightarrow\mathbb{R}:\left(
h_{1},h_{2}\right)  \mapsto b$. Hence $\Lambda_{m}\times\Lambda_{m}\ni\left(
h_{1},h_{2}\right)  \mapsto\gamma_{h_{1}^{-1}h_{2}}$ is in $L^{1}(\Lambda
_{m}\times\Lambda_{m},\mu\times\mu)$ and
\[
\lim_{n\rightarrow\infty}\int_{\Lambda_{m}\times\Lambda_{m}}A_{n}(h_{1}
,h_{2})d\left(  h_{1},h_{2}\right)
=\int_{\Lambda_{m}\times\Lambda_{m}} \gamma_{h_{1}^{-1}h_{2}}d\left(
h_{1},h_{2}\right)
\]
by Lebesgue's dominated convergence theorem. The proposition now follows by
Fubini's Theorem.
\end{proof}

Now we can finally state a van der Corput type lemma. It is worth
pointing out again that the following result continues to hold even
if $G$ is not abelian, but still amenable, and $\mu$ is right
invariant. The placing of the $g$ in (2.1.1) is then of course
important. Lemma 2.7 and Proposition 2.8 however require $\mu$ to be
left rather than right invariant when $G$ is not abelian, while they
do not directly use property (2.1.1).

\begin{theorem} 
Consider a bounded $f:G\rightarrow\mathfrak{H}$, with $\mathfrak{H}$
a Hilbert space, such that $\left\langle f(\cdot),x\right\rangle $
and $\left\langle f(\cdot ),f(\cdot)\right\rangle :G\times
G\rightarrow\mathbb{C}$ are Borel measurable (for all
$x\in\mathfrak{H}$). Assume
\[
\gamma_{h}:=\lim_{n\rightarrow\infty}\frac{1}{\mu(\Lambda_{n})}\int
_{\Lambda_{n}}\left\langle f(g),f(gh)\right\rangle dg
\]
exists for all $h\in G$. Also assume that
\begin{equation}
\lim_{n\rightarrow\infty}\frac{1}{\mu(\Lambda_{n})^{2}}\int_{\Lambda_{n}}
\int_{\Lambda_{n}}\gamma_{h_{1}^{-1}h_{2}}dh_{1}dh_{2}=0 \tag{2.6.1}
\end{equation}
(note that the integral exists by Proposition 2.5). Then
\[
\lim_{m\rightarrow\infty}\frac{1}{\mu(\Lambda_{m})}\int_{\Lambda_{m}}fd\mu=0.
\]
\end{theorem}

\begin{proof}
By Proposition 2.2 and Proposition 2.4 we only have to show that for
any $\varepsilon>0$ there is an $n$ and $m_{0}$ such that $\left|
A_{n,m}\right|  <\varepsilon$ for all $m>m_{0}$ where
\[
A_{n,m}:=\frac{1}{\mu(\Lambda_{n})^{2}}\frac{1}{\mu(\Lambda_{m})}\int
_{\Lambda_{n}}\int_{\Lambda_{n}}\int_{\Lambda_{m}}\left\langle
f(gh_{1} ),f(gh_{2})\right\rangle dgdh_{1}dh_{2}\text{.}
\]
But this follows from Proposition 2.5 and our assumptions, namely
\[
\lim_{n\rightarrow\infty}\lim_{m\rightarrow\infty}A_{n,m}=\lim_{n\rightarrow
\infty}\frac{1}{\mu(\Lambda_{n})^{2}}\int_{\Lambda_{n}}\int_{\Lambda_{n}
}\gamma_{h_{1}^{-1}h_{2}}dh_{1}dh_{2}=0.
\]
\end{proof}

We still need a few refinements regarding condition (2.6.1):

\begin{lemma} 
Let $\Lambda\subset G$ be Borel and $\mu(\Lambda)<\infty$, and
$S\subset G$ Borel such that $\Lambda^{-1} \Lambda\subset S$. For a
Borel $f:G\rightarrow\mathbb{R}^{+}$ we then have
\[
\int_{\Lambda}\int_{\Lambda}f(h_{1}^{-1}h_{2})dh_{1}dh_{2}\leq\mu(\Lambda
)\int_{S}fd\mu.
\]
\end{lemma}

\begin{proof}
Let $\chi$ denote characteristic functions, and set
$\varphi:\Lambda\times\Lambda\rightarrow G:(h_{1},h_{2})\mapsto
h_{1} ^{-1}h_{2}$. Then $f\circ\varphi$ is Borel on
$\Lambda\times\Lambda$, and therefore measurable in the product
$\sigma$-algebra on $\Lambda\times\Lambda$ obtained from $\Lambda$
's Borel $\sigma$-algebra, since $\varphi$ is continuous. Let
$Y\subset\Lambda^{-1}\Lambda$ be Borel in $G$. For $W\subset G\times
G$, let $W_{g}:=\left\{  h:(g,h)\in W\right\}  $. Then, since
$\varphi^{-1}(Y)$ is Borel in $\Lambda\times\Lambda$ and hence Borel
in $G\times G$, it follows that $\varphi^{-1}(Y)$ is in the product
$\sigma $-algebra on $G\times G$, hence we can consider
$(\mu\times\mu)\left( \varphi^{-1}(Y)\right)
=\int_{\Lambda}\mu\left(  \varphi^{-1}(Y)_{g}\right) dg $. Now
\[
\varphi^{-1}(Y)=\left\{  (g,gh):h\in Y,g\in\Lambda\cap\left(  \Lambda
h^{-1}\right)  \right\}  \subset\left\{  (g,gh):h\in Y,g\in\Lambda\right\}
=:V
\]
but $V_{g}=gY$, therefore $\mu\left(  \varphi^{-1}(Y)_{g}\right)  \leq
\mu(V_{g})=\mu(gY)=\mu(Y)$, since $\mu$ is left invariant. Hence
\begin{align*}
\int_{\Lambda\times\Lambda}\chi_{Y}\circ\varphi d(\mu\times\mu)  &
=(\mu\times\mu)\left(  \varphi^{-1}(Y)\right) \\
&  \leq\mu(\Lambda)\mu(Y)\\
&  =\mu(\Lambda)\int_{S}\chi_{Y}d\mu
\end{align*}

There is an increasing sequence $f_{n}:S\rightarrow\mathbb{R}^{+}$
of simple functions converging pointwise to $f$. From the above we
know that
\[
\int_{\Lambda\times\Lambda}f_{n}\circ\varphi d(\mu\times\mu)\leq\mu
(\Lambda)\int_{S}f_{n}d\mu
\]
and by applying Lebesgue's monotone convergence first on the right and then of
the left of this inequality, we obtain%
\[
\int_{\Lambda}\int_{\Lambda}f\left(  h_{1}^{-1}h_{2}\right)  dh_{1}dh_{2}%
=\int_{\Lambda\times\Lambda}f\circ\varphi d(\mu\times\mu)\leq\mu(\Lambda
)\int_{S}fd\mu
\]
as required, where we have used Fubini's theorem, which holds in
this case, since $f$ is non-negative.
\end{proof}

\begin{proposition} 
Consider a Borel measurable function $\gamma:G\rightarrow\mathbb{C}$
and let $\gamma_{h}:=\gamma(h)$ for all $h\in G$. Also assume that
\[
\lim_{n\rightarrow\infty}\frac{1}{\mu(\Lambda_{n})}\int_{\Lambda_{n}
^{-1}\Lambda_{n}}\left|  \gamma_{h}\right|  dh=0.
\]
Then
\[
\lim_{n\rightarrow\infty}\frac{1}{\mu(\Lambda_{n})^{2}}\int_{\Lambda_{n}}
\int_{\Lambda_{n}}\gamma_{h_{1}^{-1}h_{2}}dh_{1}dh_{2}=0
\]
if the iterated integral exists for all $n\geq n_{0}$ for some
$n_{0}$.
\end{proposition}

\begin{proof}
Since $\Lambda_{n}$ is compact, $\Lambda_{n} ^{-1}\Lambda_{n}$ is
Borel, and so
\begin{align*}
\left|
\frac{1}{\mu(\Lambda_{n})^{2}}\int_{\Lambda_{n}}\int_{\Lambda_{n}
}\gamma_{h_{1}^{-1}h_{2}}dh_{1}dh_{2}\right|   &
\leq\frac{1}{\mu(\Lambda
_{n})^{2}}\int_{\Lambda_{n}}\int_{\Lambda_{n}}\left|
\gamma_{h_{1}^{-1}h_{2}
}\right|  dh_{1}dh_{2}\\
&  \leq\frac{1}{\mu(\Lambda_{n})}\int_{\Lambda_{n}^{-1}\Lambda_{n}}\left|
\gamma_{h}\right|  dh
\end{align*}
by Lemma 2.7.
\end{proof}

\section{Weak mixing}

In this section we define weak mixing, and study some of its
characterizations using simple tools like density limits. This sets
the stage for our study of weak mixing of all orders in the next
section. As mentioned in Section 2, $G$ need not be abelian in this
section. In fact, even the properties of F\o lner sequences are not
needed until Corollary 3.8 (this means that in Definition 3.2 below
one could in principle work with an arbitrary sequence $\left(
\Lambda_{n}\right)  $ of Borel sets in $G$ with $0<\mu\left( \Lambda
_{n}\right)  <\infty$). The material in this section is fairly
standard, except that we work with the notion ``$M$-weak mixing''
(and ``$M$ -ergodicity''), which is important in Section 4.

\begin{definition} 
Let $\omega$ be a state on a unital $\ast $-algebra $A$, i.e. a
linear functional on $A$ such that $\omega(a^{\ast }a)\geq0$ and
$\omega(1)=1$. Let $\tau_{g}$ be a $\ast$-automorphism of $A$ for
every $g\in G$ such that $\tau_{g}\circ\tau_{h}=\tau_{gh}$ for all
$g,h\in G$, and such that $\tau_{e}$ is the identity on $A$ and
$G\rightarrow \mathbb{C}:g\mapsto\omega(a\tau_{g}(b))$ is Borel
measurable for all $a,b\in A$. Then we'll call $(A,\omega,\tau,G)$ a
\emph{$\ast$-dynamical system}. If we further have that $A$ is a
C*-algebra, and the state is preserved, i.e.
$\omega\circ\tau_{g}=\omega$ for all $g$ in $G$, then
$(A,\omega,\tau,G)$ is called a \emph{C*-dynamical system}.
\end{definition}

Suppose that $(A,\omega,\tau,G)$ is a $\ast$-dynamical system. Then so is
$(\overline{A},\overline{\omega},\tau,G)$ where $\overline{\omega
}(a)=\overline{\omega(a)}$, while $\overline{A}$ is the $\ast$-algebra $A$
with the original scalar multiplication replaced by $\alpha\cdot
a=\overline{\alpha}a$ for all $\alpha\in\mathbb{C}$ and $a\in A$. If
$A\otimes\overline{A}$ denotes the algebraic tensor product of $A$ with
$\overline{A}$, then $(A\otimes\overline{A},\omega\otimes\overline{\omega
},\tau\otimes\tau,G)$ is also a $\ast$-dynamical system, where $\left(
\tau\otimes\tau\right)  _{g}:=\tau_{g}\otimes\tau_{g}$. When $A$ is normed, we
assign the same norm to $\bar{A}$, and on $A\otimes\overline{A}$ for our
purposes any norm satisfying $\left\|  a\otimes b\right\|  \leq\left\|
a\right\|  \left\|  b\right\|  $ will do, for example the spatial C*-norm when
$A$ is a C*-algebra. However, even in the normed case, $A\otimes\bar{A}$ will
denote the algebraic tensor product; we will not work with the completion in
the norm.

For a group $G$, let Hom$(G)$ denote the set of \textit{all} group
homomorphisms $G\rightarrow G$.

\begin{definition} 
Let $(A,\omega,\tau,G)$ be a $\ast
$-dynamical system and consider an $M\subset$ Hom$(G)$ such that
$G\rightarrow\mathbb{C}:g\mapsto\omega(a\tau_{\varphi(g)}(b))$ is
Borel measurable for all $\varphi\in M$.

(i) $(A,\omega,\tau,G)$ is said to be $M$-\textit{weakly mixing}
\textit{relative t}o $(\Lambda_{n})$, if
\[
\lim_{n\rightarrow\infty}\frac{1}{\mu(\Lambda_{n})}\int_{\Lambda_{n}}\left|
\omega(a\tau_{\varphi(g)}(b))-\omega(a)\omega(b)\right|  dg=0
\]
for all $a,b\in A$, and for all $\varphi\in M$.

(ii) $(A,\omega,\tau,G)$ is said to be $M$-\textit{ergodic} \textit{relative
t}o $(\Lambda_{n})$, if
\[
\lim_{n\rightarrow\infty}\frac{1}{\mu(\Lambda_{n})}\int_{\Lambda_{n}}
\omega(a\tau_{\varphi(g)}(b))dg=\omega(a)\omega(b)
\]
for all $a,b\in A$, and for all $\varphi\in M$.
\end{definition}

\begin{remarks}
If $\left(  A,\omega ,\tau,G\right)  $ is a $\ast$-dynamical system,
then so is $\left( A,\omega,\tau_{\varphi(\cdot)},G\right)  $ for
any $\varphi\in$ M. So essentially we're looking at a set of systems
indexed by $M$, and one can therefore expect that the known
properties of weakly mixing and ergodic systems will extend to the
situation in Definition 3.2, as we'll review in the rest of the
section. In the case of $G=\mathbb{Z}$, $\Lambda_{n}=\left\{
1,...,n\right\}  $ and with $M=\left\{  id_{\mathbb{Z}}\right\}  $,
Definition 3.2(i) corresponds to the usual definition of weak mixing
for an action of the group $\mathbb{Z}$. Since all homomorphisms of
$\mathbb{Z}$ are of the form $n\mapsto kn$ for some
$k\in\mathbb{Z}$, one can then easily show for a C*-dynamical system
that $\left\{  id_{\mathbb{Z}}\right\}  $-weak mixing implies
Hom$(\mathbb{Z})\backslash\{0\}$-weak mixing, where $0$ is the
homomorphism $n\mapsto0$. Note more generally that if the
homomorphism given by $\varphi_{0}(g)=e$ for all $g\in G$ is in $M$
then the system cannot be expected to be $M$-weakly mixing, hence we
would not want $\varphi_{0}$ to be in $M$. We mention this simply
because $\varphi_{0}$ does appear in the theory to follow, but not
as an element of $M$.
\end{remarks}

We now turn to a few technical tools which we will need in Section 4.

\begin{definition} 
(i) A set $R\subset G$ is said to have \emph{density zero relative
to} $(\Lambda_{n})$, and we write $D_{(\Lambda _{n})}(R)=0$, if and
only if there exists a measurable set $S\subset G$, with $R\subset
S$ such that
\[
\lim_{n\rightarrow\infty}\frac{\mu(\Lambda_{n}\cap S)}{\mu(\Lambda_{n})}=0.
\]

(ii) We say that $f:G\rightarrow L$, with $L$ a real or complex
normed space, has \emph{density limit} $a\in L$ relative to
$(\Lambda_{n})$, if and only if for each $\varepsilon>0$,
$D_{(\Lambda_{n})}(S_{\varepsilon})=0$, where
\[
S_{\varepsilon}:=\{h\in G:\Vert f(h)-a\Vert\geq\varepsilon\},
\]
and we write it as
\[
D_{(\Lambda_{n})}\text{-}\lim f=D_{(\Lambda_{n})}\text{-}\lim_{h}f(h)=a.
\]
\end{definition}

Note that if $R$ and $S$ have density zero relative to $(\Lambda_{n})$ and
$V\subset S$, then $R\cap S$, $R\cup S$ and $V$ also have density zero
relative to $(\Lambda_{n})$.

We now give a Koopman-von Neumann type lemma:

\begin{lemma} 
Let $f:G\rightarrow\lbrack0,\infty)$ be bounded and measurable. Then
the following are equivalent:

(1) $D_{(\Lambda_{n})}$-$\lim f=0$

(2) $\displaystyle             \lim_{n\rightarrow\infty}\frac{1}{\mu
(\Lambda_{n})}\int_{\Lambda_{n}}f~d\mu=0$
\end{lemma}

\begin{proof}
For every $\varepsilon>0$, let $S_{\varepsilon
}:=\{h\in G:f(h)\geq\varepsilon\}$, which is a measurable set, since
$f$ is measurable.

(1) $\Rightarrow$ (2): From (1) we have that each $S_{\varepsilon}$ has
density zero relative to $(\Lambda_{n})$. Given any $\varepsilon>0$ and index
$\alpha$, consider the term
\[
\frac{1}{\mu(\Lambda_{n})}\int_{\Lambda_{n}}fd\mu=\frac{1}{\mu(\Lambda_{n}
)}\int_{\Lambda_{n}\cap
S_{\varepsilon}}f~d\mu+\frac{1}{\mu(\Lambda_{n})}
\int_{\Lambda_{n}\cap S_{\varepsilon}^{c}}fd\mu\text{.}
\]
Since $S_{\varepsilon}$ has density zero relative to $(\Lambda_{n})$
\[
0\leq\frac{1}{\mu(\Lambda_{n})}\int_{\Lambda_{n}\cap
S_{\varepsilon}} f~d\mu\leq\frac{\mu\left(  \Lambda_{n}\cap
S_{\varepsilon}\right)  } {\mu(\Lambda_{n})}\sup f(G)\rightarrow0
\]
as $n\rightarrow\infty$. Also,
\[
0\leq\frac{1}{\mu(\Lambda_{n})}\int_{\Lambda_{n}\cap
S_{\varepsilon}^{c} }f~d\mu\leq\frac{\mu\left(  \Lambda_{n}\cap
S_{\varepsilon}^{c}\right)  }
{\mu(\Lambda_{n})}\varepsilon\leq\varepsilon\newline
\]
hence
\[
\lim_{n\rightarrow\infty}\frac{1}{\mu(\Lambda_{n})}\int_{\Lambda_{n}}
fd\mu=0\text{.}
\]

(2) $\Rightarrow$ (1): Clearly $\varepsilon\chi_{S_{\varepsilon}}\leq f$. Also
note that $D_{(\Lambda_{n})}\left(  S_{\varepsilon}\right)  =0$, since
$S_{\varepsilon}$ is measurable and
\[
\varepsilon\frac{\mu\left(  \Lambda_{n}\cap S_{\varepsilon}\right)
} {\mu\left(  \Lambda_{n}\right)
}\leq\frac{1}{\mu(\Lambda_{n})}\int _{\Lambda_{n}}fd\mu
\]
which tends to zero as $n\rightarrow\infty$.
\end{proof}

\begin{corollary} 
Let $f:G\rightarrow\mathbb{R}$ be bounded and measurable. Then
\[
\lim_{n\rightarrow\infty}\frac{1}{\mu(\Lambda_{n})}\int_{\Lambda_{n}
}[f(h)]^{2}~dh=0
\]
if and only if
\[
\lim_{n\rightarrow\infty}\frac{1}{\mu(\Lambda_{n})}\int_{\Lambda_{n}
}|f(h)|~dh=0.
\]
\end{corollary}

\begin{proof}
Given any $\varepsilon>0$. Let
\[
S_{\varepsilon}:=\{h\in G:[f(h)]^{2}\geq\varepsilon^{2}\}=\{h\in
G:|f(h)|\geq\varepsilon\}.
\]
Suppose that
$\lim_{n\rightarrow\infty}\frac{1}{\mu(\Lambda_{n})}\int
_{\Lambda_{n}}[f(h)]^{2}~dh=0$, i.e. $D_{(\Lambda_{n})}$-$\lim_{h}
[f(h)]^{2}=0$ by Lemma 3.4. \noindent By the definition of the
density limit we have $D_{(\Lambda_{n})}(S_{\varepsilon})=0$. Since
$\varepsilon>0$ is arbitrary, we conclude that
$D_{(\Lambda_{n})}$-$\lim|f|=0$, and hence
$\lim_{n\rightarrow\infty}\frac{1}{\mu(\Lambda_{n})}\int_{\Lambda_{n}
}|f(h)|~dh=0$ by Lemma 3.4. The converse follows
similarly.
\end{proof}

As a result, the $\left|  \cdot\right|  $ in Definition 3.2(i) of $M$-weak
mixing, can be replaced by $\left|  \cdot\right|  ^{2}$, which is useful below
and in Section 4.

\begin{lemma} 
Let $f:G\rightarrow\mathbb{C}$ bounded and measurable. Let
$\beta\in\mathbb{C}$. If
\[
\lim_{n\rightarrow\infty}\frac{1}{\mu(\Lambda_{n})}\int_{\Lambda_{n}
}f(h)dh=\beta\mbox{ \emph{and} }  \lim_{n\rightarrow\infty}\frac{1}
{\mu(\Lambda_{n})}\int_{\Lambda_{n}}|f(h)|^{2}dh=|\beta|^{2},
\]
then
\[
\lim_{n\rightarrow\infty}\frac{1}{\mu(\Lambda_{n})}\int_{\Lambda_{n}}\left|
f(h)-\beta\right|  ^{2}dh=0.
\]
\end{lemma}

\begin{proof}
This follows immediately if we note that
\begin{align*}
&  \frac{1}{\mu(\Lambda_{n})}\int_{\Lambda_{n}}\left|  f(h)-\beta\right|
^{2}dh\\
&  =\frac{1}{\mu(\Lambda_{n})}\int_{\Lambda_{n}}( f(h)-\beta)(\overline
{f(h)-\beta})dh\\
&  =\frac{1}{\mu(\Lambda_{n})}\int_{\Lambda_{n}}\left(  |f(h)|^{2}
-\beta\overline{f(h)}-\overline{\beta}f(h)+|\beta|^{2}\right)  dh\\
&  \rightarrow0
\end{align*}
as $n\rightarrow\infty$.
\end{proof}

Next we consider standard characterizations of weak mixing, that we will need.

\begin{proposition} 
Let $(A,\omega,\tau,G)$ be a $\ast $-dynamical system and $M\subset$
Hom$(G)$ such that $g\mapsto \omega(a\tau_{\varphi(g)}(b))$ is
measurable for all $a,b\in A$ and all $\varphi\in M$. The following
are equivalent:

(1) $(A,\omega, \tau,G)$ is $M$-weakly mixing relative to
$(\Lambda_{n})$.

(2) $(A\otimes\overline{A},\omega\otimes\overline{\omega},\tau
\otimes\tau,G)$ is $M$-weakly mixing relative to $(\Lambda_{n})$.

(3) $(A\otimes\overline{A},\omega\otimes\overline{\omega},\tau
\otimes\tau,G)$ is $M$-ergodic relative to $(\Lambda_{n})$.
\end{proposition}

\begin{proof}
(2) $\Rightarrow$ (3): Follows immediately
from Definition 3.2.

(3) $\Rightarrow$ (1): Let $a,b\in A$ and $\varphi\in M$. We have
\begin{align*}
&
\lim_{n\rightarrow\infty}\frac{1}{\mu(\Lambda_{n})}\int_{\Lambda_{n}}
\omega(a\tau_{\varphi(g)}(b))dg\\
&
=\lim_{n\rightarrow\infty}\frac{1}{\mu(\Lambda_{n})}\int_{\Lambda_{n}
}\omega\otimes\overline{\omega}\left(  (a\otimes1)(\tau\otimes\tau
)_{\varphi(g)}(b\otimes1)\right)  dg\\
&  =\omega\otimes\overline{\omega}(a\otimes1)\omega\otimes\overline{\omega
}(b\otimes1)\\
&  =\omega(a)\omega(b).
\end{align*}
Furthermore,
\begin{align*}
&
\lim_{n\rightarrow\infty}\frac{1}{\mu(\Lambda_{n})}\int_{\Lambda_{n}
}\left|  \omega(a\tau_{\varphi(g)}(b))\right|  ^{2}dg\\
&
=\lim_{n\rightarrow\infty}\frac{1}{\mu(\Lambda_{n})}\int_{\Lambda_{n}
}\omega\otimes\overline{\omega}\left(  (a\otimes a)\left(
\tau\otimes
\tau\right)  _{\varphi(g)}(b\otimes b)\right)  dg\\
&  =\omega\otimes\overline{\omega}(a\otimes a)\omega\otimes\overline{\omega
}(b\otimes b)\\
&  =|\omega(a)\omega(b)|^{2}.
\end{align*}
Therefore by Lemma 3.6 we have that
\[
\lim_{n\rightarrow\infty}\frac{1}{\mu(\Lambda_{n})}\int_{\Lambda_{n}}\left|
\omega(a\tau_{\varphi(g)}(b))-\omega(a)\omega(b)\right|  ^{2}dg=0,
\]
and it follows from Corollary 3.5 that $(A,\omega,\tau,G)$ is
$M$-weakly mixing relative to $(\Lambda_{n})$.

(1) $\Rightarrow$ (2): Given any $\varphi\in M$ and $a,b\in A\otimes
\overline{A}$, with $a=\sum_{j=1}^{n}a_{j}\otimes c_{j}$ and
$b=\sum_{k=1} ^{m}b_{k}\otimes d_{k}$ where
$a_{1},\ldots,a_{n},b_{1},\ldots,b_{m}\in A$ and
$c_{1},\ldots,c_{n},d_{1},\ldots,d_{m}\in\overline{A}$, then

\begin{align*}
&  \left|  \omega\otimes\overline{\omega}\left(  a(\tau\otimes\tau
)_{\varphi(g)}(b)-\omega\otimes\overline{\omega}(a)\omega\otimes
\overline{\omega}(b)\right)  \right| \\
&  \leq\sum_{j=1}^{n}\sum_{k=1}^{m}\left(  \left|  \omega\left(  a_{j}%
\tau_{\varphi(g)}(b_{k})\right)  \overline{\omega}\left(  c_{j}\tau
_{\varphi(g)}(d_{k})\right)  -\omega\left(  a_{j}\tau_{\varphi(g)}%
(b_{k})\right)  \overline{\omega}(c_{j})\overline{\omega}(d_{k})\right|
\right. \\
&  ~~~+\left.  \left|  \omega\left(  a_{j}\tau_{\varphi(g)}(b_{k})\right)
\overline{\omega}(c_{j})\overline{\omega}(d_{k})-\omega(a_{j})\overline
{\omega}(c_{j})\omega(b_{k})\overline{\omega}(d_{k})\right|  \right) \\
&  \leq\sum_{j=1}^{n}\sum_{k=1}^{m}\left(  \Vert a_{j}\Vert_{\omega}\Vert
b_{k}\Vert_{\omega}\left|  \omega\left(  c_{j}\tau_{\varphi(g)}(d_{k})\right)
-\omega(c_{j})\omega(d_{k})\right|  \right. \\
&  ~~~+\left.  \left|  \overline{\omega}(c_{j})\overline{\omega}%
(d_{k})\right|  \left|  \omega\left(  a_{j}\tau_{\varphi(g)}(b_{k})\right)
-\omega(a_{j})\omega(b_{k})\right|  \right)
\end{align*}
hence (2) follows from Definition 3.2(i).
\end{proof}

We can now show that the definition of $M$-weak mixing relative to a F\o lner
sequence, is independent of the F\o lner sequence being used:

\begin{corollary} 
If a state preserving $\ast $-dynamical system $(A,\omega, \tau,G)$
is $M$-weakly mixing relative to some F\o lner sequence in $G$, then
it is $M$-weakly mixing relative to every F\o lner sequence in $G$.
\end{corollary}

\begin{proof}
By the mean ergodic theorem (see \cite{DDS}) the $M$-ergodicity of a
$\ast$-dynamical system is independent of the F\o lner sequence
being used. Hence $M$-weak mixing is also independent of the F\o
lner sequence by Proposition 3.7(1 and 3).
\end{proof}

\section{Weak mixing of all orders}

For a state $\omega$ on a unital $\ast$-algebra $A$, we will denote
any GNS representation of $(A,\omega)$ by $(H,\iota)$, which is a
Hilbert space $H$ and a linear mapping $\iota:A\rightarrow H$ such
that $\langle\iota (a),\iota(b)\rangle=\omega(a^{\ast}b)$ for all
$a,b\in A$ and with $\iota(A)$ dense in $H$. The corresponding
cyclic vector will be denoted by $\Omega :=\iota (1)$.  Note that we
use the convention where inner products are conjugate linear in the
first slot. For elements of $A$ we'll use the notation
$\prod_{j=1}^{k}a_{j}$ to denote the product $a_{1}...a_{k}$ in this
specific order. Lastly we remind the reader that we are still using
the group $G$ as described in Section 2, and as was mentioned there
the fact that $G$ is abelian becomes essential in the current
section.

\begin{definition} 
Let $(A,\omega,\tau,G)$ be a $\ast
$-dynamical system where $A$ has a submultiplicative norm. Such a
$\ast $-dynamical system is said to be $M$-\textit{asymptotically
abelian} \textit{relative t}o $(\Lambda_{n})$, where $M\subset$
Hom$(G)$, if $G\rightarrow A:g\mapsto\tau_{\varphi(g)}(b)$ is
continuous, and
\[
\lim_{n\rightarrow\infty}\frac{1}{\mu(\Lambda_{n})}\int_{\Lambda_{n}}
\Vert\lbrack a,\tau_{\varphi(g)}(b)]\Vert dg=0
\]
for all $a,b\in A$, and for all $\varphi\in M$, where $[a,b]:=ab-ba$.
\end{definition}

Under this assumption, we show that weak mixing implies weak mixing of all
orders. Our approach is strongly influenced by that of \cite{FKO} for the case
of a measure theoretic dynamical system and the group $\mathbb{Z}$. The proof
is by induction, two steps of which are given by the following:

\begin{proposition} 
Given $M\subset$ Hom$(G)$, let $(A,\omega,\tau,G)$ denote any
C*-dynamical system such that $G\rightarrow
A:g\mapsto\tau_{\varphi(g)}(a)$ is continuous for all $a\in A$, and
for all $\varphi\in M$. We are going to work with a collection of
such systems, but with $G$,$~\left( \Lambda_{n}\right) $ and $M$
fixed. Given $k\in\mathbb{N}$, let $\varphi_{1},\ldots ,\varphi_{k}$
denote elements of $M$ and $a_{0},\ldots,a_{k}$ elements of $A$. Set
$\varphi_{0}(h)=e$ for all $h\in G$. Consider the following
statements:

1[$k$]:
$\displaystyle\lim_{n\rightarrow\infty}\frac{1}{\mu(\Lambda_{n})}\int_{\Lambda_{n}}\left|
\omega\left(  \prod_{j=0}^{k}\tau_{\varphi_{j}(g)}(a_{j})\right)
-\prod _{j=0}^{k}\omega(a_{j})\right|  dg=0$

2[$k$]:
$\displaystyle\lim_{n\rightarrow\infty}\frac{1}{\mu(\Lambda_{n})}\int_{\Lambda_{n}}\omega\left(
\prod_{j=0}^{k}\tau_{\varphi_{j}(g)}(a_{j})\right)
dg=\prod_{j=0}^{k} \omega(a_{j})$

3[$k$]: $\displaystyle\lim_{n\rightarrow\infty}\left\|
\frac{1}{\mu(\Lambda_{n})}\int_{\Lambda_{n} }\iota\left(
\prod_{j=1}^{k}\tau_{\varphi_{j}(g)}(a_{j})\right)
dg-\prod_{j=1}^{k}\omega(a_{j})\Omega\right\|  =0$

\noindent Then

(1) 1[$k$] implies 2[$k$].

(2) If 3[$k$] holds for all $M$-weakly mixing $(A,\omega,\tau,G)$
which are $M$-asymptotically abelian relative to $(\Lambda_{n})$,
all $a_{1},...,a_{k}$ and all distinct $\varphi_{1},...,\varphi_{k}$
with $\varphi_{j}\neq\varphi_{l}$ when $j\neq l$ for $j,l\in\left\{
1,...,k\right\}  $, then 1[$k$] also holds for all $M$-weakly mixing
$(A,\omega,\tau,G)$ which are $M$-asymptotically abelian relative to
$(\Lambda_{n})$, all $a_{0},...,a_{k}$ and all distinct $\varphi
_{1},...,\varphi_{k}$ with $\varphi_{j}\neq\varphi_{l}$ when $j\neq
l$ for $j,l\in\left\{ 1,...,k\right\}  $.
\end{proposition}

\begin{proof}
(1) Trivial.

(2) Let $\kappa:=\prod_{j=1}^{k}\omega(a_{j})$. By the assumption we
have that
\begin{align*}
&
\lim_{n\rightarrow\infty}\frac{1}{\mu(\Lambda_{n})}\int_{\Lambda_{n}}
\omega\left(  \prod_{j=0}^{k}\tau_{\varphi_{j}(g)}(a_{j})\right)  ~dg\\
&  = \lim_{n\rightarrow\infty}\left\langle
\iota(a_{0}^{\ast}),\frac{1}
{\mu(\Lambda_{n})}\int_{\Lambda_{n}}\iota\left(  \prod_{j=1}^{k}\tau
_{\varphi_{j}(g)}(a_{j})~dg\right)  \right\rangle \\
&  =\left\langle \iota(a_{0}^{\ast}),\kappa\Omega\right\rangle \\
&  =\prod_{j=0}^{k}\omega(a_{j})
\end{align*}

First note that Proposition 3.7 (1 and 2) imply that the product
system
$(A\otimes\overline{A},\omega\otimes\overline{\omega},\tau\otimes\tau,G)$
is $M$-weakly mixing. Since $\left\|  a\otimes b\right\|
\leq\left\|  a\right\| \left\|  b\right\|  $, it is also easy to
verify that $(A\otimes\overline
{A},\omega\otimes\overline{\omega},\tau\otimes\tau,G)$ is
$M$-asymptotic abelian relative to $\left(  \Lambda_{n}\right)  $.
Hence by the equality above (which applies to all systems which are
both $M$-weakly mixing and $M$-asymptotic abelian) we obtain
\begin{align*}
&
\lim_{n\rightarrow\infty}\frac{1}{\mu(\Lambda_{n})}\int_{\Lambda_{n}
}\left|  \omega\left(
\prod_{j=0}^{k}\tau_{\varphi_{j}(g)}(a_{j})\right)
\right|  ^{2}dg\\
&
=\lim_{n\rightarrow\infty}\frac{1}{\mu(\Lambda_{n})}\int_{\Lambda_{n}
}\omega\otimes\overline{\omega}\left(  \prod_{j=0}^{k}(\tau\otimes
\tau)_{\varphi_{j}(g)}\left(  a_{j}\otimes a_{j}\right)  \right)  dg\\
&  =\prod_{j=0}^{k}\omega\otimes\overline{\omega}(a_{j}\otimes a_{j})\\
&  =\left|  \prod_{j=0}^{k}\omega(a_{j})\right|  ^{2}
\end{align*}
proving 1[$k$] by Corollary 3.5 and Lemma 3.6.
\end{proof}

Note that the only property of $M$-weak mixing and $M$-asymptotic abelianness
which is used in Proposition 4.2's proof, is that if a system is $M$-weakly
mixing, then so is its product system, and similarly for $M$-asymptotic
abelianness. Proposition 4.2 would still hold if we just considered $M$-weakly
mixing systems for example, or systems with some abstract property, call it
$E$, as long as the product system is again an $E$ dynamical system. $M$-weak
mixing and $M$-asymptotic abelianness will be used more directly in subsequent steps.

In order to complete the induction argument, we need 1[$1$], and
that if 2[$k-1$] holds for all relevant systems, then the same is
true for 3[$k$]. The latter requires some more work, and we will
need to specialize the $M$ that we will allow. Firstly note that
since the group $G$ is abelian, we have for any homomorphisms
$\varphi_{1}$ and $\varphi_{2}$ of $G$ that the function
$\varphi^{\prime }:G\rightarrow G$ defined by
\begin{equation}
\varphi^{\prime}(g):=\varphi_{2}(g)^{-1}\varphi_{1}(g) \tag{4.1}
\end{equation}
is also a homomorphism of $G$.

\begin{definition} 
Let $M\subset$ Hom$(G)$. We call $M$ \emph{translational} if for all
$\varphi_{1},\varphi_{2}\in M$ with $\varphi_{1}\neq\varphi_{2}$,
the homomorphism $\varphi^{\prime}$ defined by (4.1) is also in $M$.
\end{definition}

\begin{proposition} 
Let $(A,\omega,\tau,G)$ be a $C^{\ast}$-dynamical system which is
$M$-asymptotic abelian relative to $(\Lambda_{n})$, with $M$
translational. Set $\varphi_{0}(g)=e$ for all $g\in G$. Assume that
for some $k\in\mathbb{N}$
\begin{equation}
\lim_{n\rightarrow\infty}\frac{1}{\mu(\Lambda_{n})}\int_{\Lambda_{n}}
\omega\left(  \prod_{j=0}^{k-1}\tau_{\varphi_{j}(g)}(a_{j})\right)
dg=\prod_{j=0}^{k-1}\omega(a_{j}) \tag{4.4.1}
\end{equation}
for all $a_{0},...,a_{k-1}\in A$ and
$\varphi_{1},...,\varphi_{k-1}\in M$ with
$\varphi_{j}\neq\varphi_{l}$ when $j\neq l$ for $j,l\in\left\{
1,...,k-1\right\}  $, and in particular the existence of the limit
is assumed. Now set
\[
u_{h}:=\iota\left(  \prod_{j=1}^{k}\tau_{\varphi_{j}(h)}(a_{j})\right)
-\kappa\Omega
\]
for all $h\in G$, where $\kappa:=\prod_{j=1}^{k}\omega(a_{j})$, for
a given set of $a_{j}\in A$ and $\varphi_{j}\in M$ with
$\varphi_{j}\neq \varphi_{l}$ when $j\neq l$ for $j,l\in\left\{
1,...,k\right\}  $. Then
\[
\gamma_{h}:=\lim_{n\rightarrow\infty}\frac{1}{\mu(\Lambda_{n})}\int
_{\Lambda_{n}}\left\langle u_{g},u_{gh}\right\rangle dg
\]
exists (where $\left\langle \cdot,\cdot\right\rangle $ is taken in
$H$), and
\[
\gamma_{h}=\prod_{j=1}^{k}\omega\left(  a_{j}^{\ast} \tau_{\varphi
_{j}(h)}(a_{j}) \right)  -|\kappa|^{2}
\]
for all $h\in G$.
\end{proposition}

\begin{proof}
We have
\begin{multline*}
\left\langle u_{g},u_{gh}\right\rangle =\omega\left(  \left(
\prod_{j=1} ^{k}\tau_{\varphi_{j}(g)}(a_{j})\right)  ^{\ast}\left(
\prod_{j=1}^{k}
\tau_{\varphi_{j}(gh)}(a_{j})\right)  \right) \\
-\kappa\overline{\omega\left(  \prod_{j=1}^{k}\tau_{\varphi_{j}(g)}
(a_{j})\right)  }-\overline{\kappa}\omega\left(  \prod_{j=1}^{k}\tau
_{\varphi_{j}(gh)}(a_{j})\right)  +|\kappa|^{2}
\end{multline*}

We now consider each of the terms in the last expression separately:

(a) We see that
\[
\omega\left(  \left(  \prod_{j=1}^{k}\tau_{\varphi_{j}(g)}(a_{j})\right)
^{\ast}\left(  \prod_{j=1}^{k}\tau_{\varphi_{j}(gh)}(a_{j})\right)  \right)
=\omega\left(  \prod_{j=-k}^{k}\tau_{\varphi_{|j|}(g)}(b_{j})\right)
\]
where $b_{j}=\left\{
\begin{array}
[c]{ll}
a_{|j|}^{\ast}, & \hbox{if $-k\leq j\leq -1$;} \\
1, & \hbox{if $j=0$;} \\
\tau_{\varphi_{j}(h)}(a_{j}), & \hbox{if $1\leq j\leq k$.}
\end{array}
\right.  $

For $T_{-k},...,T_{k}\in A$ with $\left\|  T_{j}\right\|  \leq c$, one has
(see Lemma 7.4 in \cite{NSZ}) that
\[
\left\|  \prod_{j=-k}^{k}T_{j}-T_{0}\left(  \prod_{j=1}^{k}\left(
T_{-j} T_{j}\right)  \right)  \right\|  \leq
c^{2k-1}\sum_{r=1}^{k}\sum _{l=-(r-1)}^{r-1}\left\|  \left[
T_{-r},T_{l}\right]  \right\|
\]
and applying this to $T_{j}=\tau_{\varphi_{\left|  j\right|
}(g)}\left( b_{j}\right)  $ with $c:=\max_{-k\leq j\leq k}\left\|
b_{j}\right\|  $, and keeping in mind that $\left\|
\tau_{g}(a)\right\|  =\left\|  a\right\|  $, it follows that
\[
\left\|  \prod_{j=-k}^{k}\tau_{\varphi_{\left|  j\right| }(g)}\left(
b_{j}\right)  -\prod_{j=1}^{k} \tau_{\varphi_{j}(g)}\left(  b_{-j}
b_{j}\right)  \right  \|  \leq c^{2k-1}\sum_{r=1}^{k}\sum
_{l=-(r-1)}^{r-1}\left\|  \left[  b_{-r},\tau_{\varphi_{r}(g)^{-1}
\varphi_{\left|  l\right|  }(g)}\left(  b_{l}\right)  \right]
\right\|
\]
Since $M$ is translational, and $(A,\omega,\tau,G)$ is asymptotically abelian
relative to $\left(  \Lambda_{n}\right)  $, it follows that
\begin{align*}
&
\lim_{n\rightarrow\infty}\frac{1}{\mu(\Lambda_{n})}\int_{\Lambda_{n}}
\omega\left(  \prod_{j=-k}^{k}\tau_{\varphi_{|j|}(g)}(b_{j})\right)  dg\\
&
=\lim_{n\rightarrow\infty}\frac{1}{\mu(\Lambda_{n})}\int_{\Lambda_{n}
}\omega\left(  \prod_{j=1}^{k}\tau_{\varphi_{j}(g)}(b_{-j}b_{j})\right)  dg\\
&
=\lim_{n\rightarrow\infty}\frac{1}{\mu(\Lambda_{n})}\int_{\Lambda_{n}
}\omega\left(  \tau_{\varphi_{1}(g)}\left(
\prod_{j=1}^{k}\tau_{\varphi
_{1}(g)^{-1}\varphi_{j}(g)}(b_{-j}b_{j})\right)  \right)  dg\\
&
=\lim_{n\rightarrow\infty}\frac{1}{\mu(\Lambda_{n})}\int_{\Lambda_{n}
}\omega\left(
\prod_{j=0}^{k-1}\tau_{\varphi_{j}^{\prime}(g)}(b_{-(j+1)}
b_{j+1})\right)  dg\\
&  =\prod_{j=0}^{k-1}\omega\left(  b_{-(j+1)}b_{j+1}\right) \\
&  =\prod_{j=1}^{k}\omega\left(  a_{j}^{\ast}\tau_{\varphi_{j}(h)}
(a_{j})\right)
\end{align*}
by using (4.4.1), where
$\varphi_{j}^{\prime}(g):=\varphi_{1}(g)^{-1} \varphi_{j+1}(g)$ for
all $g\in G$ and $j=0,\ldots,k-1$, so $\varphi _{j}^{\prime}\in M$
for $j=1,\ldots,k-1$ since $M$ is translational. Note that
$\varphi_{j}^{\prime}\neq\varphi_{l}^{\prime}$ when $j\neq l$ for
$j,l\in\{1,\ldots,k-1\}$ and $\varphi_{0}^{\prime}(g)=e$ for all
$g\in G$, as required in our assumption. \noindent Hence
\[
\lim_{n\rightarrow\infty}\frac{1}{\mu(\Lambda_{n})}\int_{\Lambda_{n}}
\omega\left(  \left(
\prod_{j=1}^{k}\tau_{\varphi_{j}(g)}(a_{j})\right) ^{\ast}\left(
\prod_{j=1}^{k}\tau_{\varphi_{j}(gh)}(a_{j})\right)  \right)
dg=\prod_{j=1}^{k}\omega\left(  a_{j}^{\ast}\tau_{\varphi_{j}(h)}
(a_{j})\right)
\]

(b) It follows as in (a) that
\begin{align*}
\lim_{n\rightarrow\infty}\frac{1}{\mu(\Lambda_{n})}\int_{\Lambda_{n}}
\overline{\omega\left(
\prod_{j=1}^{k}\tau_{\varphi_{j}(g)}(a_{j})\right)  } dg  &
=\overline{\lim_{n\rightarrow\infty}\frac{1}{\mu(\Lambda_{n})}
\int_{\Lambda_{n}}\omega\left(
\prod_{j=0}^{k-1}\tau_{\varphi_{j}^{\prime
}(g)}(a_{j+1})\right)  dg}\\
&  =\overline{\prod_{j=0}^{k-1}\omega(a_{j+1})}\\
&  =\overline{\kappa}%
\end{align*}
by assumption.

(c) Lastly, using similar arguments as before,
\begin{align*}
&
\lim_{n\rightarrow\infty}\frac{1}{\mu(\Lambda_{n})}\int_{\Lambda_{n}
}\omega\left(  \prod_{j=1}^{k}\tau_{\varphi_{j}(gh)}(a_{j})\right)  dg\\
&
=\lim_{n\rightarrow\infty}\frac{1}{\mu(\Lambda_{n})}\int_{\Lambda_{n}
}\omega\left(  \prod_{j=0}^{k-1}\tau_{\varphi_{j}^{\prime}(g)}\left(
\tau_{\varphi_{j+1}(h)}(a_{j+1})\right)  \right)  dg\\
&  =\prod_{j=0}^{k-1}\omega(\tau_{\varphi_{j+1}(h)}(a_{j+1}))\\
&  =\prod_{j=0}^{k-1}\omega(a_{j+1})\\
&  =\kappa.
\end{align*}

(d) From (a)-(c)
\[
\gamma_{h}=\prod_{j=1}^{k}\omega\left(  a^{\ast}_{j}\left(
\tau_{\varphi _{j}(h)}(a_{j})\right)  \right)  -|\kappa|^{2}
\]
and in particular $\gamma_{h}$ exists.
\end{proof}

Next we prove that weak mixing implies weak mixing of all orders. This is
where our van der Corput lemma is finally applied, along with Propositions 4.2
and 4.4. We need one more definition:

\begin{definition} 
A F\o lner sequence $(\Lambda_{n})$ in $G$ is
said to satisfy the \emph{Tempelman condition} if there is a real
number $c>0$ such that
\[
\mu(\Lambda_{n}^{-1}\Lambda_{n})\leq c\mu(\Lambda_{n})
\]
for all $n\in \mathbb{N}$.
\end{definition}

See \cite{L} for some discussion and further references related to this condition.

\begin{theorem} 
Assume that there exists a F\o lner sequence $(\Lambda_{n})$ in $G$,
satisfying the Tempelman condition, and such that
$(\Lambda_{n}^{-1}\Lambda_{n})$ is also a F\o lner sequence in $G$.
Let $M\subset$ Hom$(G)$ be translational. Let $(A,\omega ,\tau,G)$
be an $M$-weakly mixing $C^{\ast}$-dynamical system which is
$M$-asymptotically abelian relative to $(\Lambda_{n})$. Assume
furthermore that $G\rightarrow A:g\mapsto\tau_{\varphi(g)}(a)$ is
continuous in the norm topology on $A$ for all $\varphi\in M$ and
all $a\in A$. Then
\[
\lim_{n\rightarrow\infty}\frac{1}{\mu(\Lambda_{n})}\int_{\Lambda_{n}}\left|
\omega\left(  \prod_{j=0}^{k}\tau_{\varphi_{j}(g)}(a_{j})\right)  -\prod
_{j=0}^{k}\omega(a_{j})\right|  dg=0
\]
for any $a_{j}\in A$ and any $\varphi_{1},...,\varphi_{k}\in M$ with
$\varphi_{j}\neq\varphi_{l}$ when $j\neq l$ for $j,l\in\left\{
1,...,k\right\}  $, and with $\varphi_{0}(g)=e$ for all $g\in G$.
\end{theorem}

\begin{proof}
We need to complete the induction argument started in
Proposition 4.2, and we will continue using its notation and that of
Proposition 4.4. Since $G\rightarrow A:g\mapsto\tau_{\varphi(g)}(a)$
is continuous, so is $G\rightarrow
A:g\mapsto\prod_{j=1}^{k}\tau_{\varphi_{j}(g)}(a_{j})$ in the
norm-topology on $A$. We then also have that $G\rightarrow
H:g\mapsto\iota\left(  \prod
_{j=1}^{k}\tau_{\varphi_{j}(g)}(a_{j})\right)  $ is continuous,
since $\left\|  \iota(a)\right\|  \leq\left\|  a\right\|  $. It
follows that
\[
G\times G\rightarrow\mathbb{R}:(g,h)\mapsto\left\langle \iota\left(
\prod_{j=1}^{k}\tau_{\varphi_{j}(g)}(a_{j})\right)  ,\iota\left(
\prod _{j=1}^{k}\tau_{\varphi_{j}(h)}(a_{j})\right)  \right\rangle
\]
is continuous and hence $G\times G\rightarrow\mathbb{C}:(g,h)\mapsto
\left\langle u_{g},u_{h}\right\rangle $ is continuous and therefore
Borel measurable. Note that $g\mapsto\left\langle
u_{g},x\right\rangle $ is also Borel measurable for all $x\in H$.
Furthermore, $G\rightarrow H:g\mapsto u_{g}$ is bounded. (We need
these properties, since we will be applying Theorem 2.6 to the
function $g\mapsto u_{g}$.) Since $\mu(\Lambda_{n}
^{-1}\Lambda_{n})\leq c\mu(\Lambda_{n})$, and we have $M$-weak
mixing relative to $\left(  \Lambda_{n}^{-1}\Lambda_{n}\right)  $ by
Corollary 3.8, it follows that
\begin{equation}
\lim_{n\rightarrow\infty}\frac{1}{\mu(\Lambda_{n})}\int_{\Lambda_{n}
^{-1}\Lambda_{n}}\left|  \omega\left(  a\tau_{\varphi(g)}(b)\right)
-\omega(a)\omega(b)\right|  dg=0 \tag{4.6.1}
\end{equation}
for all $a,b\in A$ and $\varphi\in M$. By Proposition 4.4, assuming
2[$k-1$] for all $M$-weakly mixing C*-dynamical systems, which are
$M$-asymptotically abelian relative to $\left(  \Lambda_{n}\right)
$ (for the given $G$ and $\left(  \Lambda_{n}\right)  $), and of
course for all $a_{0},...,a_{k-1}$ and all
$\varphi_{1},...,\varphi_{k-1}\in M$ with
$\varphi_{j}\neq\varphi_{l}$ when $j\neq l$ for $j,l\in\left\{
1,...,k-1\right\}  $, we have
\[
\gamma_{h}:=\lim_{n\rightarrow\infty}\frac{1}{\mu(\Lambda_{n})}\int
_{\Lambda_{n}}\left\langle u_{g},u_{gh}\right\rangle
dg=\prod_{j=1}^{k} \omega\left(  a_{j}^{\ast}\left(
\tau_{\varphi_{j}(h)}(a_{j})\right) \right)
-\prod_{j=1}^{k}|\omega(a_{j})|^{2}
\]
for any $a_{1},...,a_{k}$ and all $\varphi_{1},...,\varphi_{k}\in M$ with
$\varphi_{j}\neq\varphi_{l}$ when $j\neq l$ for $j,l\in\left\{
1,...,k\right\}  $, for all $h\in G$. Using the following identity (also see
Section 4 in \cite{FKO})
\begin{equation}
\prod_{j=1}^{k}c_{j}-\prod_{j=1}^{k}d_{j}=\sum_{j=1}^{k}\left(  \prod
_{l=1}^{j-1}c_{l}\right)  \left(  c_{j}-d_{j}\right)  \left(  \prod
_{l=j+1}^{k}d_{l}\right)  \tag{4.6.2}%
\end{equation}
which holds in any algebra and is easily verified by induction, it follows
that
\[
\int_{\Lambda_{m}^{-1}\Lambda_{m}}\left|  \gamma_{h}\right|  dh\leq\sum
_{j=1}^{k}A_{j}\prod_{l=j+1}^{k}|\omega(a_{l})|^{2}\int_{\Lambda_{m}%
^{-1}\Lambda_{m}}\left|  \omega\left(  a_{j}^{\ast}\left(  \tau_{\varphi
_{j}(h)}(a_{j})\right)  \right)  -|\omega(a_{j})|^{2}\right|  dh
\]
where $A_{j}:=\sup_{h\in G}\left|  \prod_{l=1}^{j-1}\omega\left(  a_{l}^{\ast
}(\tau_{\varphi_{l}(h)}(a_{l})\right)  \right|  \leq\prod_{l=1}^{j-1}\Vert
a_{l}\Vert^{2}$.

Note that $\int_{\Lambda_{m}^{-1}\Lambda_{m}}\left|  \gamma_{h}\right|  dh$
exists, since the integrand is continuous. Hence
\[
\lim_{m\rightarrow\infty}\frac{1}{\mu(\Lambda_{m})}\int_{\Lambda_{m}
^{-1}\Lambda_{m}}\left|  \gamma_{h}\right|  dh=0
\]
by (4.6.1). From Proposition 2.8 and Theorem 2.6 we then have
\[
\lim_{n\rightarrow\infty}\frac{1}{\mu(\Lambda_{n)}}\int_{\Lambda_{n}}
u_{g}dg=0
\]
i.e., 3[$k$] holds for all $M$-weakly mixing C*-dynamical systems,
which are $M$-asymptotically abelian relative to $(\Lambda_{n})$,
and all $a_{1},...,a_{k}$ and all $\varphi _{1},...,\varphi_{k}\in
M$ with $\varphi_{j}\neq\varphi_{l}$ when $j\neq l$
for $j,l\in\left\{  1,...,k\right\}  $. But 1[$1$] holds for all $a_{0}%
,a_{1}\in A$ and all $\varphi\in M$ for all $M$-weakly mixing C*-dynamical
systems, which are $M$-asymptotically abelian relative to $(\Lambda_{n})$, by
Definition 3.2(i), completing the induction argument started in Proposition
4.2, and proving 1[$k$] for all $k\in\mathbb{N}$.
\end{proof}

We now briefly consider simple examples of F\o lner sequences with the
required properties.

In the case where $G=\mathbb{Z}$ with the counting measure $\mu$, and
$\Lambda_{n}=\{-n,\ldots,n\}$ which is F\o lner in $\mathbb{Z}$, we have
$\Lambda_{n}^{-1}\Lambda_{n}=\{-2n,\ldots,2n\}$ which is also F\o lner, and
$\mu(\Lambda_{n})\leq\mu(\Lambda_{n}^{-1}\Lambda_{n})\leq2\mu(\Lambda_{n})$
for $n\geq1$. Similarly for $\mathbb{Z}^{q}$.

As another example, let $\Lambda_{m}$ be the closed ball of radius $m$ in
$\mathbb{R}^{q}$ for any positive integer $q$. Note that $\left(  \Lambda
_{m}\right)  $ is a F\o lner sequence in $\mathbb{R}^{q}$ with $\Lambda
_{m}^{-1}\Lambda_{m}=\Lambda_{2m}$, and $\mu(\Lambda_{m}^{-1}\Lambda
_{m})=2^{q}\mu(\Lambda_{m})$.

Concerning the assumption that $M$ is translational, a simple
example would be of the following type: Use the group
$G=\mathbb{R}^{q}$. Let $M$ be all $q\times q$ non-zero diagonal
real matrices acting as linear operators on $\mathbb{R}^{q}$. (We
exclude the zero matrix simply because this would make $M$-weak
mixing impossible.) Then $M$ is a translational set of homomorphisms
of $\mathbb{R}^{q}$. The same is true if we drop the condition that
the matrices be diagonal. Similarly if we work with $\mathbb{Z}^{q}$
instead of $\mathbb{R}^{q}$ and use matrices over the integers.
These examples work simply because if $\varphi_1,\varphi_2\in M$ and
$\varphi_1\neq\varphi_2$ then
$-\varphi_2(g)+\varphi_1(g)=(\varphi_1-\varphi_2)(g)$ while
$\varphi_1-\varphi_2\in M$.

\section{Compact systems}

In this section we prove a Szemer\'{e}di type property for compact
C*-systems as defined in Definition 5.1 below. Again we follow the
basic structure of the proof given in \cite{FKO}, but we have to
take into account certain subtleties and technical difficulties
arising from working with a noncommutative C*-algebra rather than
with the abelian algebra $L^{\infty}(\nu)$ used in \cite{FKO}, and
with more general groups and semigroups than $\mathbb{Z}$ and
$\mathbb{N}$. First some notation and terminology.

A linear functional $\omega$ on a $\ast$-algebra $A$ is called positive if
$\omega(a^{\ast}a)\geq0$ for all $a\in A$. This allows us to define a seminorm
$\Vert\cdot\Vert_{\omega}$ on $A$ by
\[
\Vert a\Vert_{\omega}:=\sqrt{\omega(a^{\ast}a)}%
\]
for all $a\in A$, as is easily verified using the Cauchy-Schwarz inequality
for positive linear functionals.

A set $V$ in a pseudo metric space $(X,d)$ is said to be $\varepsilon
$\emph{-separated}, where $\varepsilon>0$, if $d(x,y)\geq\varepsilon$ for all
$x,y\in V$ with $x\neq y$. A set $B\subset X$ is said to be \emph{totally
bounded in} $\left(  X,d\right)  $ if for every $\varepsilon>0$ there exists a
finite set $M_{\varepsilon}\subset X$, called a finite $\varepsilon
$\textit{-net}, such that for every $x\in B$ there is a $y\in M_{\varepsilon}$
with $d(x,y)<\varepsilon$. It is then not difficult to show that for any
$\varepsilon>0$ there exists a maximal set (in the sense of cardinality, or
number of elements) $V\subset B$ that is $\varepsilon$-separated, and
furthermore, if $B\neq\varnothing$, then $V$ is finite with $|V|>0$.

\begin{definition} 
Let $\omega$ be a positive linear functional
on a $\ast$-algebra $A$, $K$ a semigroup, and $\tau_{g}:A\rightarrow
A$ a linear map for each $g\in K$ such that
\[
\tau_{g}\circ\tau_{h}=\tau_{gh}
\]
and
\[
\Vert\tau_{g}(a)\Vert_{\omega}=\Vert a\Vert_{\omega}
\]
for all $g,h\in K$ and $a\in A$. Assume that the \emph{orbit}
\[
B_{a}:=\{\tau_{g}(a):g\in K\}
\]
is totally bounded in $(A,\Vert\cdot\Vert_{\omega})$ for each $a\in
A$. Then we call $\left(  A,\omega,\tau,K\right)  $ a \emph{compact
system}. If furthermore $A$ is a C*-algebra and
$\Vert\tau_{g}(a)\Vert\leq\Vert a\Vert$ in $A$'s norm for all $a\in
A$ and $g\in K$, then we refer to $\left( A,\omega,\tau,K\right)  $
as a \emph{compact C*-system}.
\end{definition}

In particular, if the orbits of the $\ast$-dynamical systems and C*-dynamical
systems in Definition 3.1 are totally bounded in $(A,\Vert\cdot\Vert_{\omega
})$, then those systems will be called \textit{compact}. Using the GNS
construction, it is not too difficult to see that the $L^{2}$ definition of
compactness mentioned in Section 1 is a special case of Definition 5.1.

In this paper weakly mixing systems and compact systems appear as
part of the structure of ergodic systems, so concrete examples of
weakly mixing systems and compact systems are not crucial for our
goal. Nevertheless, it is interesting to look at an example of a
compact C*-dynamical system in which the C*-algebra is
noncommutative. To do this we need a few simple tools, which we now
discuss.

First note that if a set in a C*-algebra $A$ is totally bounded in
$A$ (i.e. in terms of $A$'s norm), then it is also totally bounded
in $\left(  A,\left| \left|  \cdot\right|  \right|  _{\omega}\right)
$ for any positive linear functional $\omega$ on $A$, since $\left|
\left|  \cdot\right|  \right| _{\omega}\leq\left|  \left|
\omega\right|  \right|  ^{1/2}\left|  \left| \cdot\right|  \right| $
(keep in mind that $\omega$ is bounded, since it is positive and $A$
is a C*-algebra). Hence, if we can prove that the orbits of a given
C*-dynamical system $\left(  A,\omega,\tau,K\right)  $ are totally
bounded in $A$, then it follows that the system is compact. Of
course, this is then a stronger form of compactness, but Example 5.4
happens to possess this stronger property, and it turns out to be
easier to prove this than to prove compactness directly in terms of
$\left|  \left|  \cdot\right|  \right| _{\omega}$, since $A$'s norm
is submultiplicative, which makes it easier to work with than
$\left|  \left|  \cdot\right|  \right|  _{\omega}$.

In Lemma 5.2 and Proposition 5.3 below, we work with a C*-algebra
$A$, an arbitrary \textit{set} $K$, and a $\ast$-homomorphism
$\tau_{g}:A\rightarrow A$ for each $g\in K$. When we say that an
``orbit'' $(a_{g})\equiv\left( a_{g}\right)  _{g\in K}$ is
\textit{totally bounded} in a space, we mean that the set
$\{a_{g}:g\in K\}$ is totally bounded in that space. For any subset
$\mathfrak{V}\subset A$ we will denote the set of all polynomials
over $\mathbb{C}$ generated by the elements of $\mathfrak{V}$ and
their adjoints, by $p(\mathfrak{V})$, i.e. $p(\mathfrak{A})$
consists of all finite linear
combinations of all finite products of elements of $\mathfrak{V}%
\cup\mathfrak{V}^{\ast}$ with $\mathfrak{V}^{\ast}:=\left\{  a^{\ast}%
:a\in\mathfrak{V}\right\}  $. We will use the notation $XY:=\left\{  xy:x\in
X,y\in Y\right\}  $ whenever $X$ and $Y$ are sets for which this
multiplication of their elements is defined.

\begin{lemma} 
If $\left(  \tau_{g}(a)\right) $ is totally bounded in $A$ for every
$a$ in some subset $\mathfrak{V}$ of $A$, then $\left(
\tau_{g}(a)\right) $ is totally bounded in $A$ for every $a\in
p(\mathfrak{V})$.
\end{lemma}

\begin{proof}
Consider any $a,b\in A$ for which $\left(  \tau
_{g}(a)\right)  $ and $\left(  \tau_{g}(b)\right)  $ are totally
bounded in $A$, and any $\varepsilon>0$. By the hypothesis there are
finite sets $M,N\subset A$ such that for each $g\in K$ there is an
$a_{g}\in M$ and a $b_{g}\in N$ such that
$\Vert\tau_{g}(a)-a_{g}\Vert<\varepsilon$ and
$\Vert\tau_{g}(b)-b_{g}\Vert<\varepsilon$. Clearly
\begin{align*}
\Vert\tau_{g}(a)\tau_{g}(b)-a_{g}b_{g}\Vert &  \leq\Vert\tau_{g}(a)\Vert
\Vert\tau_{g}(b)-b_{g}\Vert+\Vert\tau_{g}(a)-a_{g}\Vert\Vert b_{g}\Vert\\
&  \leq\varepsilon\left(  \Vert\tau_{g}(a)\Vert+\Vert b_{g}\Vert\right)
\end{align*}
but note that $\Vert\tau_{g}(a)\Vert\leq\Vert a\Vert$, since $\tau_{g}$ is a
$\ast$-homomorphism and $A$ is a C*-algebra, while $\left|  \left|
b_{g}\right|  \right|  <\Vert\tau_{g}(b)\Vert+\varepsilon\leq\left|  \left|
b\right|  \right|  +\varepsilon$. Since $MN$ is a finite subset of $A$, and
$a_{g}b_{g}\in MN$, it follows that $(\tau_{g}(ab))$ is totally bounded in
$A$. Similarly $(\tau_{g}(a^{\ast}))$ and $(\tau_{g}(\alpha a+\beta b))$ are
totally bounded in $A$ for any $\alpha,\beta\in\mathbb{C}$, and this is enough
to prove the lemma.
\end{proof}

\begin{proposition} 
Now assume that $A$ is generated by a subset $\mathfrak{V}\subset A$
for which $\tau _{g}(\mathfrak{V})\subset p(\mathfrak{V})$ for every
$g\in K$. Also assume that $(\tau_{g}(a))$ is totally bounded in $A$
for every $a\in\mathfrak{V}$. Then $(\tau_{g} (a))$ is totally
bounded in $A$ for every $a\in A$.
\end{proposition}

\begin{proof}
Firstly it is easily shown that if $Y$ is a dense
subspace of a normed space $X$, $U_{g}:Y\rightarrow Y$ is linear
with $\Vert U_{g}\Vert\leq1$ for all $g\in K$, and $\left(
U_{g}y\right)  $ is totally bounded in $X$ for every $y\in Y$ (or in
$Y$ for every $y\in Y$), then for the unique bounded linear
extension $U_{g}:X\rightarrow X$ the ``orbit'' $\left( U_{g}x\right)
$ is totally bounded in $X$ for every $x\in X$. (We also used this
fact when we discussed the GNS-construction above.)

Now simply set $X=A$, $Y=p(\mathfrak{V})$ and $U_{g}=\tau_{g}$, then by our
assumptions and Lemma 5.2 all the requirements in the remark above are met.
\end{proof}

\begin{example} 
We consider a so-called rotation C*-algebra, and use Proposition 5.3
to show that we obtain a compact C*-dynamical system. As described
in Chapter VI of \cite{D}, let $\mathfrak{H}:=L^{2}
(\mathbb{R}/\mathbb{Z})$ and define two unitary operators $U$ and
$V$ on $\mathfrak{H}$ by
\[
\left(  Uf\right)  (t)=f(t+\theta)
\]
and
\[
\left(  Vf\right)  (t)=e^{2\pi it}f(t)
\]
for $f\in\mathfrak{H}$, where $\theta\in\mathbb{R}$ (though the
interesting case is $\theta\notin\mathbb{Q}$). These operators
satisfy
\begin{equation}
UV=e^{2\pi i\theta}VU\text{.} \tag{5.4.1}%
\end{equation}
Let $\mathfrak{A}$ be the C*-algebra generated by $U$ and $V$. Note that
$\mathfrak{A}$ is noncommutative because of (5.4.1). Then, as shown in Chapter
VI of \cite{D}, there is a unique trace $\omega$ on $\mathfrak{A}$, i.e. a
state with $\omega(ab)=\omega(ba)$. Define $\tau:\mathfrak{A}\rightarrow
\mathfrak{A}$ by $\tau(a)=U^{\ast}aU$ for all $a\in\mathfrak{A}$, then $\tau$
is a $\ast$-isomorphism and therefore $\left|  \left|  \tau(a)\right|
\right|  =\left|  \left|  a\right|  \right|  $, since $\mathfrak{A}$ is a
C*-algebra. Also, since $\omega$ is a trace and $U$ is unitary, $\Vert
\tau(a)\Vert_{\omega}=\Vert a\Vert_{\omega}$ for all $a\in\mathfrak{A}$. Hence
$(\mathfrak{A},\omega,\tau,\mathbb{N})$ is a C*-dynamical system, where by
slight abuse of notation $\tau$ here denotes the function $n\mapsto\tau^{n}$
as well, to fit it into Definitions 3.1 and 5.1's notation.

We now show that $(\mathfrak{A},\omega,\tau,\mathbb{N})$ is compact: It is
trivial that $(\tau^{n}(U))=(U)$ is totally bounded in $\mathfrak{A}$.
Furthermore, $\tau^{n}(V)=(U^{\ast})^{n}VU^{n}=e^{-2\pi in\theta}V$ by
(5.4.1). Since the unit circle is compact, it follows that $(\tau^{n}(V))$ is
totally bounded in $\mathfrak{A}$. From Proposition 5.3 with $\mathfrak{V}%
=\left\{  U,V\right\}  $ we conclude that $(\tau^{n}(a))$ is totally bounded
in $\mathfrak{A}$ for all $a\in\mathfrak{A}$. In particular $(\mathfrak{A}%
,\omega,\tau,\mathbb{N})$ is a compact C*-system. Similarly $(\mathfrak{A}%
,\omega,\tau,\mathbb{Z})$ is a compact C*-dynamical system.
\end{example}

Now we resume the general theory.

\begin{definition} 
Let $K$ be a semigroup. We call a set
$E\subset K$ \emph{relatively dense in }$K$ if there exist an
$r\in\mathbb{N}$ and $g_{1},\ldots,g_{r}\in K$ such that
\[
E\cap\{gg_{1},\ldots,gg_{r}\}\neq\varnothing
\]
for all $g\in K$.
\end{definition}

Strictly speaking one could call this \textit{left} relative denseness, with
the right hand case being defined similarly in terms of $g_{j}g$, but we will
only work with Definition 5.5 in this paper. The usual definition of relative
denseness of a subset $E$ in $\mathbb{N}$ is in terms of ``bounded gaps'' (see
\cite{P} for example), and it is easy to check that in this special case the
two definitions are equivalent.

\begin{proposition} 
Let $K$ be a semigroup, $\left(  X,\Vert\cdot\Vert\right)  $ a
seminormed space, and $U_{g}:X\rightarrow X$ a linear map for each
$g\in K$ such that $U_{g}U_{h}=U_{gh}$ and $\Vert
U_{g}x\Vert\geq\Vert x\Vert $ for all $g,h\in K$ and $x\in X$.
Suppose that $B_{x_{0}}:=\{U_{g}x_{0}:g\in K\}$ is totally bounded
in $\left( X,\Vert\cdot\Vert\right)  $ for some $x_{0}\in X$. Then
for each $\varepsilon>0$, the set
\[
E:=\{g\in K:\Vert U_{g}x_{0}-x_{0}\Vert<\varepsilon\}
\]
is relatively dense in $K$.
\end{proposition}

\begin{proof}
Since $B_{x_{0}}$ is totally bounded in $\left(
X,\Vert\cdot\Vert\right)  $, there is a maximal $V=\{U_{g_{1}}x_{0}
,...,U_{g_{r}}x_{0}\}$, with $U_{g_{j}}x_{0}\neq U_{g_{l}}x_{0}$
whenever $j\neq l$, which is $\varepsilon$-separated. \noindent But
$\Vert
U_{g^{\prime}gg_{j}}x_{0}-U_{g^{\prime}gg_{l}}x_{0}\Vert\geq\Vert
U_{g_{j} }x_{0}-U_{g_{l}}x_{0}\Vert$ for any $g,g^{\prime}\in K$,
hence $V_{g^{\prime
}g}:=\{U_{g^{\prime}gg_{1}}x_{0},...,U_{g^{\prime}gg_{r}}x_{0}\}$ is
$\varepsilon$-separated, with $r$ elements. Since
$V_{g^{\prime}g}\subset B_{x_{0}}$, it is also maximally
$\varepsilon$-separated in $B_{x_{0}}$. But $U_{g^{\prime}}x_{0}\in
B_{x_{0}}$, therefore $\Vert U_{gg_{j}}x_{0} -x_{0}\Vert\leq\Vert
U_{g^{\prime}gg_{j}}x_{0}-U_{g^{\prime}}x_{0} \Vert<\varepsilon$ for
some $j\in\{1,\ldots,r\}$. The last inequality follows from the
maximality of $V_{g^{\prime}g}$ in the following way. Suppose this
inequality wasn't true for some $j$, then
$V_{g^{\prime}g}\cup\left\{ U_{g^{\prime}}x_{0}\right\}  $ would be
$\varepsilon$-separated and contained in $B_{x_{0}}$, but with
strictly greater cardinality than $V_{g^{\prime}g}$, since it would
have one more element, namely $U_{g^{\prime}}x_{0}$ (which of course
is not already in $V_{g^{\prime}g}$ if the inequality doesn't hold
for any $j$), contradicting maximality.

Hence, for each $g\in K$ there exists an
$h\in\{gg_{1},\ldots,gg_{r}\}$ such that $\Vert U_{h}
x_{0}-x_{0}\Vert<\varepsilon$, i.e.
\[
E\cap\{gg_{1},\ldots,gg_{r}\}\neq\varnothing
\]
for all $g\in K$, and so $E$ is relatively dense in $K$.
\end{proof}

\begin{corollary} 
Let $\left(  A,\omega,\tau,K\right) $ be a compact system and let
$m_{0},...,m_{k}\in\mathbb{N}\cup \{0\}$. For any $\varepsilon>0$
and $a\in A$, the set
\[
E:=\{g\in K:\Vert\tau_{g^{m_{j}}}(a)-a\Vert_{\omega}<\varepsilon\text{ for
}j=0,...,k\}
\]
is then relatively dense in $K$, where we write
$\tau_{g^{0}}(a)\equiv a$.
\end{corollary}

\begin{proof}
Without loss we can assume that none of the $m_{j}$ 's
are zero. Then the result follows from Proposition 5.6 with
$\varepsilon$ replaced by $\varepsilon/\max\{m_{0},\ldots,m_{k}\}$,
since for every $j=0,...,k$ we have
\begin{align*}
&  \Vert\tau_{g^{m_{j}}}(a)-a\Vert_{\omega}\\
&  \leq\Vert\tau_{g^{m_{j}}}(a)-\tau_{g^{m_{j}-1}}(a)\Vert_{\omega}+\Vert
\tau_{g^{m_{j}-1}}(a)-\tau_{g^{m_{j}-2}}(a)\Vert_{\omega}+\ldots+\Vert\tau
_{g}(a)-a\Vert_{\omega}\\
&  =\Vert\tau_{g^{m_{j}-1}}[\tau_{g}(a)-a]\Vert_{\omega}+\Vert\tau
_{g^{m_{j}-2}}[\tau_{g}(a)-a]\Vert_{\omega}+\ldots+\Vert\tau_{g}
(a)-a\Vert_{\omega}\\
&  =m_{j}\Vert\tau_{g}(a)-a\Vert_{\omega}\\
&  <\varepsilon
\end{align*}
for all $g\in K$ for which $\Vert\tau_{g}(a)-a\Vert_{\omega}<\varepsilon
/\max\{m_{0},\ldots,m_{k}\}$.
\end{proof}

A positive linear functional $\omega$ on a C*-algebra $A$ is bounded, and
without loss we can assume that $\left|  \left|  \omega\right|  \right|  =1$
(the case $\omega=0$ being trivial), i.e. $\omega$ is a \textit{state} on $A$.
By the Cauchy-Schwarz inequality we have
\[
\left|  \omega(ab)\right|  \leq\left|  \left|  a^{\ast}\right|  \right|
_{\omega}\left|  \left|  b\right|  \right|  _{\omega}\leq\sqrt{\left|  \left|
aa^{\ast}\right|  \right|  }\left|  \left|  b\right|  \right|  _{\omega
}=\left|  \left|  a\right|  \right|  \left|  \left|  b\right|  \right|
_{\omega}%
\]
A \textit{trace} is defined to be a state $\omega$ on a C*-algebra
$A$ such that $\omega(ab)=\omega(ba)$ for all $a,b\in A$. Note that
from the previous inequality we then we have
\[
|\omega(abc)|=|\omega(cab)|\leq\Vert a\Vert\Vert b\Vert_{\omega}\Vert c\Vert
\]
for all $a,b,c\in A$. This fact is used in the proof of Lemma 5.8. The set of
positive elements of $A$ will be denoted by $A^{+}$.

\begin{lemma} 
Let $A$ be a C*-algebra and $\omega$ a trace on $A$. Suppose that
$b\in A^{+}$, $\Vert b\Vert\leq1$ and $\omega(b)>0$. Let $k\in
\mathbb{N}\cup\{0\}$, then $\omega(b^{k+1})>0$ so we can choose
$\varepsilon>0$ such that $\varepsilon<\omega(b^{k+1})$. Consider
$c_{0},\ldots,c_{k}\in A$ such that $\Vert c_{j}\Vert\leq1$ and
$\Vert c_{j}-b\Vert_{\omega}<\varepsilon /(k+1)$ for $j=0,\ldots,k$.
Then
\[
\left|  \omega\left(  \prod_{j=0}^{k}c_{j}\right)  \right|
>\omega\left( b^{k+1}\right)  -\varepsilon>0\text{.}
\]
\end{lemma}

\begin{proof}
We clearly have $\omega(b^{k+1})>0$. Furthermore,
\begin{align*}
\left|  \omega\left(  \prod_{j=0}^{k}c_{j}\right)  -\omega(b^{k+1})\right|
&  =\left|  \omega\left(  \prod_{j=0}^{k}c_{j}-\prod_{j=0}^{k}b\right)
\right| \\
&  \leq\sum_{j=0}^{k}\left(  \left\|  \prod_{l=0}^{j-1}c_{l}\right\|  \left\|
c_{j}-b\right\|  _{\omega}\left\|  b^{k-j}\right\|  \right) \\
&  <\varepsilon.
\end{align*}
where we've used (4.6.2).
\end{proof}

\begin{corollary} 
Let $(A,\omega,\tau,K)$ be a compact C*-system with $\omega$ a
trace. Suppose that $a\in A^{+},$ and $\omega(a)>0$. Take any
$m_{0},\ldots,m_{k}\in\mathbb{N}\cup\{0\}$ and any $\varepsilon>0$
with $\varepsilon<\omega\left( a^{k+1}\right)$. Then there exists a
relatively dense set $E$ in $K$ such that
\[
\left|  \omega\left(  \prod_{j=0}^{k}\tau_{g^{m_{j}}}(a)\right)  \right|
>\omega\left(  a^{k+1}\right)  -\varepsilon>0
\]
for all $g\in E$.
\end{corollary}

\begin{proof}
Since $\omega(a)>0$, $\left|  \left|  a\right| \right|
>0$, so we can set $b:=a/\left|  \left|  a\right|  \right|  $. For
$c_{j}:=\tau_{g^{m_{j}}}(b)$ we have $\left|  \left|  c_{j}\right|
\right| \leq\left|  \left|  b\right|  \right|  =1$, so from Lemma
5.8 it follows that
\[
\left|  \omega\left(  \prod_{j=0}^{k}\tau_{g^{m_{j}}}(b)\right)  \right|
>\omega\left(  b^{k+1}\right)  -\frac{\varepsilon}{\left\|  a\right\|
^{k+1}}
\]
for every $g\in K$ for which $\Vert\tau_{g^{m_{j}}}(b)-b\Vert_{\omega
}<\varepsilon/\left\|  a\right\|  ^{k+1}(k+1)$ for all $j=0,...,k$. By
Corollary 5.7 this set of $g$ 's is relatively dense in $K$.
\end{proof}

So far in this section we have not used F\o lner sequences. However,
for the remainder of this section we again make use of such
sequences, so let $G$ and $\left(  \Lambda_{n}\right)  $ be as in
Section 2; in particular $G$ is abelian. We can make a simple
refinement without complicating the proofs, namely in the rest of
this section let $K$ be a Borel set in $G$ which forms a
(necessarily abelian) semigroup and contains each $\Lambda_{n}$. We
then say that $\left( \Lambda_{n}\right)  $ is a F\o lner sequence
in $K$. (A trivial example is $G=\mathbb{Z}$, $\Lambda_{n}
=\{1,...,n\}$ and $K=\mathbb{N}$.) This is just to make clear that
only semigroup structure is used in this section. In the next
section, where Theorem 5.13 below is applied, we will of course take
$K$ to be $G$. Let $\Sigma$ denote the $\sigma$-algebra of Borel
sets of $G$ that are contained in $K$.

\begin{remark}
It is also interesting to note that the arguments below do not
require the semigroup $K$ to be abelian, however this would require
the existence of a slightly different type of F\o lner sequence.
Namely, if $G$ were non-abelian, and $\mu$ right invariant, one
would have to assume the existence of a sequence (or a net) of Borel
sets $\left( \Lambda_{n}\right)  $ of $G$ (which are contained in
$K$) with $0<\mu (\Lambda_{n})<\infty$ such that
$\lim_{n\rightarrow\infty}\mu\left(
\Lambda_{n}\Delta(g\Lambda_{n})\right)  /\mu(\Lambda_{n})=0$ or all
$g\in K$. Note that $g$ is to the left of the set despite $\mu$
being right invariant.
\end{remark}

Before we reach the main result of this section, we have to discuss F\o lner
sequences a bit further.

\begin{lemma} 
Take any $g_{n}\in K$ for each $n$. Then the sequence
\[
\left(  \Lambda_{n}g_{n}\right)
\]
is also a F\o lner sequence in $K$.
\end{lemma}

\begin{proof}
Since $K$ has the right cancellation property, we have
$(Ag)\Delta(Bg)=(A\Delta B)g$ for all $A,B\subset K$ and $g\in K$.\
Hence
\begin{align*}
\frac{\mu\left(  (\Lambda_{n}g_{n})\Delta(g(\Lambda_{n}g_{\alpha}))\right)
}{\mu(\Lambda_{n}g_{n})}  &  =\frac{\mu\left(  (\Lambda_{n}\Delta(g\Lambda
_{n}))g_{n}\right)  }{\mu(\Lambda_{n}g_{n})}\\
&  =\frac{\mu\left(  \Lambda_{n}\Delta(g\Lambda_{n})\right)  }{\mu(\Lambda
_{n})}\\
&  \longrightarrow0
\end{align*}
as $n\rightarrow\infty$.
\end{proof}

\begin{definition} 
Let $\left(  \Lambda_{n}\right)  $ be any
F\o lner sequence in $K$. Consider any $V\in\Sigma$ and set
\[
D_{\left(  \Lambda_{n}\right)  }(V):=\lim_{n\rightarrow\infty}\left[
\inf\left\{  \frac{\mu(\Lambda_{m}\cap V)}{\mu(\Lambda_{m})}:m\geq
n\right\} \right]
\equiv\liminf_{n\rightarrow\infty}\frac{\mu(\Lambda_{n}\cap
V)}{\mu(\Lambda_{n})}.
\]
If $D_{\left(  \Lambda_{n}\right)  }(V)>0$, then we say that\ $V$
has \emph{positive lower density relative to
}$\left(\Lambda_{n}\right) $.
\end{definition}

It is easily checked that $D_{\left(  \Lambda_{n}\right)  }(V)$ in
this definition always exists.

\begin{lemma} 
Let $E\in\Sigma$ be relatively dense in $K$. Then:

(1) There exists an $r\in\mathbb{N}$ and $g_{1},\ldots ,g_{r}\in K$
such that the following holds: for each $B\in\Sigma $ with
$\mu(B)<\infty$ there exists a $j\in\{1,\ldots ,r\}$ such that
$\mu((Bg_{j})\cap E)\geq\mu(B)/r$.

(2) $E$ has positive lower density relative to some F\o lner
sequence in $K$.

(3) Let $f:K\rightarrow\mathbb{R}$ a $\Sigma$-measurable function
with $f\geq0$. Assume that $f(g)\geq\alpha$ for some $\alpha>0$ and
all $g\in E\in\Sigma$. Then there exists a F\o lner sequence
$(\Lambda_{n} )$ in $K$ such that
\[
\liminf_{n\rightarrow\infty}\frac{1}{\mu(\Lambda_{n})}\int_{\Lambda_{n}}
fd\mu>0\text{.}
\]
\end{lemma}

\begin{proof}
(1) Let $g_{1},...,g_{r}$ be given by Definition 5.5. Set
$B_{j}:=\{b\in B:bg_{j}\in E\}$ for $j=1,\ldots,r$, so $B_{j}
g_{j}=(Bg_{j})\cap E\in\Sigma$ and hence $B_{j}\in\Sigma$. Now, for
any $b\in B$ we know from Definition 2.3 that
$E\cap\{bg_{1},\ldots,bg_{r} \}\neq\varnothing$. So $bg_{j}\in E$
for some $j\in\{1,\ldots,r\}$, i.e. $b\in B_{j}$. Hence
$B=\bigcup_{j=1}^{r}B_{j}$ and therefore
\[
\mu(B)=\mu(\bigcup_{j=1}^{r}B_{j})\leq\sum_{j=1}^{r}\mu(B_{j})=\sum_{j=1}
^{r}\mu(B_{j}g_{j})=\sum_{j=1}^{r}\mu((Bg_{j})\cap E)
\]
from which the conclusion follows.

(2) Consider any F\o lner sequence $\left(  \Lambda_{n}\right)  $ in $K$. Let
$g_{1},\ldots,g_{r}\in K$ be as in Definition 5.5. \noindent For each $n$ it
follows from (1) that there exists a $j(n)\in\{1,\ldots,r\}$ such that
\[
\frac{\mu((\Lambda_{n}g_{j(n)})\cap
E)}{\mu(\Lambda_{n}g_{j(n)})}\geq \frac{1}{r}
\]
where we also made use of $\mu(\Lambda_{n}g_{j(n)})=\mu(\Lambda_{n})$. But it
follows from Lemma 5.10 that $\left(  \Lambda_{n}^{\prime}\right)  $ given by
$\Lambda_{n}^{\prime}:=\Lambda_{n}g_{j(n)}$ is a F\o lner sequence in $K$.
Furthermore,
\[
D_{\left(  \Lambda_{n}^{\prime}\right)  }(E)=\liminf_{n\rightarrow\infty
}\frac{\mu(\Lambda_{n}^{\prime}\cap E)}{\mu(\Lambda_{n}^{\prime})}\geq
\lim_{n\rightarrow\infty}\frac{1}{r}=\frac{1}{r}.
\]

(3) By (2) there exists a F\o lner sequence $(\Lambda_{n})$ in $K$ such that
\[
\liminf_{n\rightarrow\infty}\frac{1}{\mu(\Lambda_{n})}\int_{\Lambda_{n}}
fd\mu\geq\liminf_{n\rightarrow\infty}\frac{1}{\mu(\Lambda_{n})}\int
_{\Lambda_{n}\cap E}\alpha~dg=\alpha D_{\left(  \Lambda_{n}\right)
}(E)>0.
\]
\end{proof}

Finally we reach the goal of this section, namely a Szemer\'{e}di
type property for compact $C^{\ast}$-systems:

\begin{theorem} 
Let $(A,\omega,\tau,K)$ be a compact C*-system with $\omega$ a trace
and $K$ a Borel measurable semigroup in $G$ such that
$\Lambda_{n}\subset K$ for every $n$, where $G$ and
$\left(\Lambda_{n}\right) $ are as in Section 2. Let $a\in A^{+}$
with $\omega(a)>0$. Take any $m_{0},\ldots,m_{k}\in
\mathbb{N}\cup\{0\}$. Assume that $g\mapsto\omega\left(
\prod_{j=0}^{k}\tau_{g^{m_{j}}}(a)\right) $ and
$g\mapsto\Vert\tau_{g^{m_{j}}}(a)-a\Vert_{\omega}$ are
$\Sigma$-measurable on $K$ for $j=0,1,\ldots,k$. Then there exists a
F\o lner sequence $(\Lambda_{n}^{\prime})$ in $K$ such that
\[
\liminf_{n\rightarrow\infty}\frac{1}{\mu(\Lambda_{n}^{\prime})}\int
_{\Lambda_{n}^{\prime}}\left|  \omega\left(  \prod_{j=0}^{k}\tau_{g^{m_{j}}%
}(a)\right)  \right|  d\mu(g)>0.
\]
\end{theorem}

\begin{proof}
This follows from Lemma 5.12(3) and Corollary 5.9,
since $E=\{g\in
K:\Vert\tau_{g^{m_{j}}}(a)-a\Vert_{\omega}<\varepsilon\text{ for
}j=0,...,k\}$ is $\Sigma$-measurable.
\end{proof}

Note that if for example we assume that $g\mapsto\tau_{g}(a)$ is continuous in
$A$'s norm, then both $g\mapsto\omega\left(  \prod_{j=0}^{k}\tau_{g^{m_{j}}%
}(a)\right)  $ and $g\mapsto\Vert\tau_{g^{m_{i}}}(a)-a\Vert_{\omega}$ are
continuous and hence Borel measurable.

\section{Ergodic systems}

In measure theoretic ergodic theory it is well known that a system
is weakly mixing if and only if it contains no non-trivial compact
factors. In this section we show that this result can be extended to
noncommutative ergodic theory, where the measure space (and its
algebra of $L^{\infty}$-functions) is replaced by a $\sigma$-finite
von Neumann algebra and a faithful normal state. To avoid confusion
we stress that the word ``factor'' as used in this paper does not
refer to a von Neumann algebra which is a factor (i.e. has trivial
center), but to a subsystem of a dynamical system as defined below.
The two main ingredients of the proof are a so-called ``proper value
theorem'' due to St\o rmer (Theorem 2.5 in \cite{S}), and the
splitting theorem of Jacobs-Deleeuw-Glicksberg (see Section 2.4 in
\cite{K}). Once this is done, we prove our final result regarding
the Szemer\'{e}di property in ergodic systems.

In this section $G$ is a completely arbitrary group for which we
will state additional requirements (like being abelian or locally
compact) as needed.

\begin{definition} 
A \emph{W*-dynamical system} $\left( A,\omega,\tau,G\right)  $
consists of a von Neumann algebra $A$ on which we have a faithful
normal state $\omega$, and where $\tau:G\rightarrow
$Aut$(A):g\mapsto\tau_{g}$ is a representation of any abelian group
$G$ as *-automorphisms of $A$ (i.e. $\tau_{e}=$ id$_{A}$ and
$\tau_{g}\circ\tau _{h}=\tau_{gh}$), such that
$\omega\circ\tau_{g}=\omega$ for all $g$ in $G$.
\end{definition}

Note that the existence of a faithful normal state on $A$ in this definition
implies that $A$ is $\sigma$-finite (see for example Proposition 2.5.6 in
\cite{BR}). It is convenient to work in the GNS representation of such a
system and for certain intermediate results the group $G$ need not be abelian,
therefore we will mostly work with the following:

\begin{definition} 
A \emph{represented system} $\left( R,\omega_{\Omega},\alpha\right)$
consists of the following: Firstly a von Neumann algebra $R$ on a
Hilbert space $H$, a unit vector $\Omega\in H$ which is cyclic and
separating for $R$, in terms of which we define a state
$\omega_{\Omega}$ on $R$ by $\omega_{\Omega}(a)=\left\langle \Omega
,a\Omega\right\rangle $. Furthermore we have a unitary
representation $U:G\rightarrow B(H):g\mapsto U_{g}$ of an arbitrary
group $G$ (i.e. $U_{g}$ is a unitary operator, $U_{e}=1$ and
$U_{g}U_{h}=U_{gh}$), such that $U_{g}\Omega=\Omega$ and
$U_{g}RU_{g}^{\ast}\subset R$ for all $g\in G$, and in terms of
which $\alpha:G\rightarrow$Aut$(R):g\mapsto\alpha_{g}$ is defined by
$\alpha_{g}(a)=U_{g}aU_{g}^{\ast}$.
\end{definition}

The notation in these two definitions will be used consistently, for
example reference to a represented system will imply the notation
$\left( R,\omega_{\Omega},\alpha\right)  $, $H$, $G$ and $U$, and
throughout the rest of this section $\left(
R,\omega_{\Omega},\alpha\right)  $ is a represented system. Note
that the GNS representation $\left(  H,\pi,\Omega\right)  $ of a
W*-dynamical system $\left(  A,\omega,\tau,G\right)  $ gives us a
\textit{corresponding} represented system $\left(  R,\omega_{\Omega}
,\alpha\right)  $ where $R=\pi(A)$ and
$\alpha_{g}(\pi(a))=\pi(\tau_{g}(a))$ in terms of which $U$ is
uniquely defined. Also keep in mind that $\pi$ is faithful in this
situation.

For a represented system an \textit{eigenoperator} of $\alpha$ is an
$a\in A\backslash\{0\}$ such that there exists a function
$\lambda_{a} :G\rightarrow\mathbb{C}$ with
$\alpha_{g}(a)=\lambda_{a}(g)a$ for all $g\in G$. Note that in this
case $\left|  \lambda_{a}(g)\right| =1$ for all $g\in G$, i.e.
$\lambda_{a}$ is  \textit{unimodular} and hence a group homomorphism
to the circle. Similarly an \textit{eigenvector} of $U$ is an $x\in
H\backslash\{0\}$ such that there exists a function
$\lambda_{x}:G\rightarrow\mathbb{C}$ with $U_{g}x=\lambda_{x}(g)x$
for all $g$. Again note that $\lambda_{x}$ is unimodular and hence a
group homomorphism to the circle.

For a represented system we will denote the Hilbert subspace of $H$ spanned by
the eigenvectors of $U$ by $H_{0}$. The Hilbert subspace of $H$ spanned by the
eigenvectors $x$ with $\lambda_{x}=1$ will be denoted by $H_{1}$. Note that
$\mathbb{C}\Omega\subset H_{1}\subset H_{0}$, with equality allowed.

\begin{definition} 
A represented system is called \emph{ergodic} (respectively
\emph{weakly mixing}) when $\dim H_{1}=1$ (respectively $\dim
H_{0}=1$). A W*-dynamical system is called \emph{ergodic}
(respectively \emph{weakly mixing}) when its corresponding
represented system is ergodic (respectively weakly mixing).
\end{definition}

When $G$ is as in Section 2, then ergodicity and weak mixing of the dynamical
system $\left(  A,\omega,\alpha\right)  $ as given in Definition 6.3 are
equivalent to $\left\{  \text{id}_{G}\right\}  $-ergodicity and $\left\{
\text{id}_{G}\right\}  $-weak mixing as given in Definition 3.2. For
ergodicity this follows from the mean ergodic theorem, and for weak mixing it
can be shown to follow from the general theory in Section 2.4 of \cite{K}.

For a represented system we define a norm $\left\|  \cdot\right\|  _{\Omega}$
on $R$ by $\left\|  a\right\|  _{\Omega}=\omega_{\Omega}(a^{\ast}%
a)^{1/2}=\left\|  a\Omega\right\|  $ where $\left\|  \cdot\right\|  $ denotes
the norm of $H$.

\begin{definition} 
A \emph{factor} $\left( N,\omega,\tau\right) $ of a W*-dynamical
system $\left(  A,\omega ,\tau,G\right)  $ consists of a
$\ast$-algebra $N\subset A$ and the restrictions of $\omega$ and
$\tau_{g}$ to $N$, such that $\tau_{g}(N)\subset N$ for all $g$.
Similarly a \emph{factor} $\left(  N,\omega_{\Omega} ,\alpha\right)
$ of a represented system $\left(  R,\omega_{\Omega} ,\alpha\right)
$ consists of a $\ast$-algebra $N\subset R$ and the restrictions of
$\omega_{\Omega}$ and $\alpha_{g}$ to $N$, such that
$\alpha_{g}(N)\subset N$ for all $g$. Such factors are called
\emph{compact} if respectively every orbit $\tau_{G}(a)=\left\{
\tau_{g}(a):g\in G\right\} $ is totally bounded in $\left( N,\left\|
\cdot\right\|  _{\omega}\right)  $ or every orbit
$\alpha_{G}(a)=\left\{  \alpha_{g}(a):g\in G\right\}  $ is totally
bounded in $\left(  N,\left\|  \cdot\right\|  _{\Omega}\right)  $. A
factor will be called \emph{nontrivial} if $R$ strictly contains
$\mathbb{C}1$.
\end{definition}

For a represented system the orbit of $x\in H$ will be denoted by
$U_{G}x=\left\{  U_{g}x:g\in G\right\}  $.

Let $B(H)$ denote the algebra of all bounded linear operators in the Hilbert
space $H$, and let $S^{\prime}$ denote the commutant of a set $S\subset B(H)$.

Let $B$ denote the set of all eigenoperators of $\alpha$ together with the
zero operator, and let $C$ be the $\ast$-algebra generated by $B$. Set
$N=C^{\prime\prime}$, hence $N$ is a von Neumann algebra contained in $R$.

\begin{proposition} 
$\left(  N,\omega_{\Omega},\alpha\right) $ is a compact factor of
$\left( R,\omega_{\Omega},\alpha\right)  $.
\end{proposition}

\begin{proof}
It is easily seen that $B$ is closed under adjoints,
products and scalar multiples, hence
\begin{equation}
C=\left\{  \sum_{j=1}^{n}a_{j}:a_{1},...,a_{n}\in B,n>0\right\}
\tag{6.5.1}
\end{equation}
but for $a\in B$ we have $\alpha_{h}(\alpha_{g}(a))=\lambda_{a}(h)\alpha
_{g}(a)$, hence $\alpha_{g}(B)\subset B$ and $\alpha_{g}(C)\subset C$ for all
$g$. For any $a\in N$ there exists a net $(a_{\gamma})$ in $C$ such that
$a_{\gamma}x\rightarrow ax$ for all $x\in H$ according to von Neumann's
density theorem. Therefore $\left\langle x,\alpha_{g}(a_{\gamma}%
)y\right\rangle =\left\langle U_{g}^{\ast}x,a_{\gamma}U_{g}^{\ast
}y\right\rangle \rightarrow\left\langle U_{g}^{\ast}x,aU_{g}^{\ast
}y\right\rangle =\left\langle x,\alpha_{g}(a)y\right\rangle $ for
all $x,y\in H$ and all $g\in G$. In other words
$\alpha_{g}(a_{\gamma})$ converges in the weak operator topology to
$\alpha_{g}(a)$, but $\alpha_{g}(a_{\gamma})\in C\subset N$ so by
the bicommutant theorem $\alpha_{g}(a)\in N$. This proves that
$\alpha_{g}(N)\subset N$, and therefore $\left(  N,\omega_{\Omega}
,\alpha\right)  $ is a factor of $\left(
R,\omega_{\Omega},\alpha\right)  $.

Next we show that $\left(  N,\omega_{\Omega},\alpha\right)  $ is
compact. First note that for $a\in B$ we have
$\alpha_{G}(a)=\lambda_{a}(G)a$ or $\alpha_{G}(a)=\{0\}$, and both
these orbits are totally bounded in $B$ with the pseudo metric
obtained by restricting $\left\|  \cdot\right\|  _{\Omega}$ to $B$,
since $\lambda_{a}(G)$ is a subset of the unit circle (which is
compact) in $\mathbb{C}$. From (6.5.1) we can then conclude that
$\alpha _{G}(a)$ is totally bounded in $\left(  C,\left\|
\cdot\right\|  _{\Omega }\right)  $ for every $a\in C$. Again by von
Neumann's density theorem for any $a\in N$ and $\varepsilon>0$ there
is $b\in C$ such that $\left\| a-b\right\|  _{\Omega}=\left\|
a\Omega-b\Omega\right\|  <\varepsilon$. Now, if $E$ is a finite
$\varepsilon$-net in $\left(  C,\left\|  \cdot\right\|
_{\Omega}\right)  $ for $\tau_{G}(b)$, then from
\[
\left\|  \alpha_{g}(a)-c\right\|  _{\Omega}\leq\left\|  \alpha_{g}%
(a)-\alpha_{g}(b)\right\|  _{\Omega}+\left\|  \alpha_{g}(b)-c\right\|
_{\Omega}<2\varepsilon
\]
for some $c\in E$ we see that $E$ is a finite $2\varepsilon$-net for
$\alpha_{G}(a)$, i.e. the latter is totally bounded.
\end{proof}

\begin{lemma} 
If $G$ is abelian, then $H_{0}$ is the set of all elements $x\in H$
whose orbits $U_{G} x$ in $H$ are totally bounded.
\end{lemma}

\begin{proof}
This is essentially a special case of general results
proven in Section 2.4 of \cite{K}. We show how it follows from those
general results. Let $H_{k}$ be the set of all elements of $H$ with
totally bounded orbits under $U$.

Let $E$ be the set of all eigenvectors of $U$, so $H_{0}=\overline
{\text{span}E}$. Clearly $U_{G}x$ is totally bounded for every $x\in E$, since
$\lambda_{x}(G)$ is a subset of the unit circle in $\mathbb{C}$. From this it
follows that $U_{G}x$ is totally bounded for every $x\in H_{0}$. I.e.
$H_{0}\subset H_{k}.$

Now suppose that $H_{0}\neq H_{k}$. It is straightforward to show that $H_{k}$
is a Hilbert subspace of $H$, so it follows that there an $x\in H_{k}%
\backslash\{0\}$ which is orthogonal to $H_{0}$. So $x\in H_{v}:=H\ominus
H_{0}$, but $H_{v}$ is the space of so-called ``flight vectors'' which means
that there is an $S$ in the weak operator closure $\overline{U_{G}}$ of
$U_{G}$ in $B(H)$ such that $Sx=0$. This is the essence of the splitting
theorem in Section 2.4 of \cite{K} as applied to a unitary group on a Hilbert
space. However, it is easily seen that $U_{G}y$ is relatively weakly compact
for every $y\in H$, since any closed ball in $H$ is weakly compact, hence
according to Lemma 2.4.2 in \cite{K} $\overline{U_{G}}x=\overline{U_{G}x}^{w}$
where $\overline{U_{G}x}^{w}$ denotes the weak closure of $U_{G}x$ in $H$. The
norm closure $\overline{U_{G}x}$ of the totally bounded set $U_{G}x$ is
compact and therefore weakly compact and hence weakly closed, so
$\overline{U_{G}x}^{w}\subset\overline{U_{G}x}$. Putting all this together we
have $0=Sx\in\overline{U_{G}x}$ contradicting $x\neq0$ and the fact that $U$
is a unitary group and hence normpreserving.
\end{proof}

\begin{proposition} 
(1) If $\left(  R,\omega_{\Omega },\alpha\right)  $ is ergodic but
not weakly mixing, then the factor $\left(
N,\omega_{\Omega},\alpha\right)  $ is nontrivial.

(2) If $\left(  R,\omega_{\Omega},\alpha\right)  $ is weakly mixing
and $G$ is abelian, then every element $a\in M$ with a totally
bounded orbit $\alpha_{G}(a)$ lies in $\mathbb{C} 1$. In particular
$\left(  R,\omega_{\Omega},\alpha\right)  $ has no nontrivial
compact factor.
\end{proposition}

\begin{proof}
(1) According to Theorem 2.5 in \cite{S} the map
$a\mapsto a\Omega$ is a bijection from the set of eigenoperators of
$\alpha$ to the set of eigenvectors of $U$, where we simultaneously
note that our definition of ergodicity of $\left(
R,\omega_{\Omega},\alpha\right)  $ is equivalent to that of \cite{S}
(namely $\alpha_{g}(a)=a$ for all $g$ implies that
$a\in\mathbb{C}1)$, since $\Omega$ is cyclic and separating for $R$.
Since $H_{0}$ strictly contains $\mathbb{C}\Omega$ by Definition
6.3, it follows that $B$ strictly contains $\mathbb{C}1$, hence $N$
strictly contains $\mathbb{C}1$. Therefore $\left(
N,\omega_{\Omega},\alpha\right)  $ is indeed nontrivial.

(2) Consider any $a\in R$ with totally bounded orbit, then from
$\left\| \alpha_{g}(a)-b\right\|  _{\Omega}=\left\|
U_{g}a\Omega-b\Omega\right\|  $ we see that $U_{G}a\Omega$ is
totally bounded in $H$. So $a\Omega\in \mathbb{C}\Omega$ by Lemma
6.6 and Definition 6.3, but $\Omega$ is separating for $R$, hence
$a\in\mathbb{C}1$. In particular the $\ast$-algebra of any compact
factor of $\left(  R,\omega_{\Omega},\alpha\right)  $ must be
contained in $\mathbb{C}1$.
\end{proof}

Using these results, we can now prove

\begin{theorem} 
Let $\left(  A,\omega,\tau,G\right) $ be an ergodic W*-dynamical
system. Then $\left(  A,\omega ,\tau,G\right) $ is weakly mixing if
and only if it has no non-trivial compact factor.
\end{theorem}

\begin{proof}
Let $\left(  H,\pi,\Omega\right)  $ be the GNS
representation of $\left(  A,\omega\right)  $ and $\left(
R,\omega_{\Omega },\alpha\right)  $ the corresponding represented
system. Since $\pi$ is faithful (giving a $\ast$-isomorphism
$A\rightarrow R$), it is simple to show that a nontrivial compact
factor in $\left(  A,\omega,\tau,G\right)  $ gives one in $\left(
R,\omega_{\Omega},\alpha\right)  $, and vice versa. The theorem then
follows from Propositions 6.5 and 6.7.
\end{proof}

\begin{definition} 
A W*-dynamical system $\left( A,\omega
,\tau,G\right)  $, with $G$ as in Section 2, is said to have the
\emph{Szemer\'{e}di property} if there exists a F\o lner sequence
$(\Lambda_{n})$ in $G$ such that for any $k\in\mathbb{N}$ and
$m_{1},\dots,m_{k}\in\mathbb{N}$ with $m_{1}<\ldots <m_{k}$ and for
all $a\in A^{+}$ with $\omega(a)>0$,
\[
\liminf_{n\rightarrow\infty}\frac{1}{\mu(\Lambda_{n})}\int_{\Lambda_{n}%
}\left|  \omega\left(  a\prod_{j=1}^{k}\tau_{g^{m_{j}}}(a)\right)
\right| dg>0.
\]
\end{definition}

\begin{theorem} 
Suppose that $\left(A,\omega,\tau,G\right) $ is a non-trivial (i.e.
$A\neq\mathbb{C}$) ergodic W*-dynamical system, with $G$ locally
compact, second countable, (and abelian, by Definition 6.1), and
containing a F\o lner sequence $(\Lambda_{n})$ satisfying the
Tempelman condition and such that $(\Lambda_{n}^{-1}\Lambda_{n})$ is
also a F\o lner sequence. Let $\omega$ be a trace and $g\mapsto
\tau_{g}(a)$ be continuous for every $a\in A$. Suppose that for each
$m\in\mathbb{Z}\backslash\{0\}$ there exists a F\o lner sequence
$(\Gamma_{n})$ and a $c>0$ such that
\[
\frac{1}{\mu(\Lambda_{n})}\int_{\Lambda_{n}}f(g^{m})~dg\leq\frac{c}
{\mu(\Gamma_{n})}\int_{\Gamma_{n}}f(g)dg
\]
for all Borel measurable $f:G\rightarrow [0,\infty)$ and all $n$,
and such that $\left( A,\omega,\tau,G\right) $ is
$\{id_G\}$-asymptotically abelian relative to $(\Gamma_{n})$. Then
the Szemer\'{e}di property holds for a nontrivial factor of $\left(
A,\omega,\tau,G\right)  $.
\end{theorem}

\begin{proof}
Set $m_0 :=0$ and $g^0 :=e$. Let $a\in A$ with $\omega(a)>0$. From
Theorem 6.8 it follows that $\left(A,\omega,\tau,G\right)  $ is
either weakly mixing or it has a non-trivial compact factor. If
$\left(  A,\omega,\tau,G\right)  $ is weakly mixing then
$\left(A,\omega,\tau,G\right)$ is $M$-weakly mixing, where
$M:=\{\varphi_{m}:m\in\mathbb{Z}\backslash\{0\}\}$ and $\varphi
_{m}:G\rightarrow G:g\mapsto g^{m}$. This can be seen from
\[
\frac{1}{\mu(\Lambda_{n})}\int_{\Lambda_{n}}\left|
\omega(a\tau_{g^{m} }(b))-\omega(a)\omega(b)\right|dg
\leq\frac{c}{\mu(\Gamma_{n})} \int_{\Gamma_{n}}\left|
\omega(a\tau_{g}(b))-\omega(a)\omega (b)\right|dg  \rightarrow0
\]
as $n\rightarrow\infty$. Similarly $\left(A,\omega,\tau,G\right)$ is
$M$-asymptotically abelian relative to $(\Lambda_{n})$. Since $M$ is
translational, it now follows from Theorem 4.6 that
\[
\lim_{n\rightarrow\infty}\frac{1}{\mu(\Lambda_{n})}\int_{\Lambda_{n}}\left|
\omega\left(  \prod_{j=0}^{k}\tau_{g^{m_{j}}}(a)\right)  -\omega
(a)^{k+1}\right|  dg=0.
\]
and therefore
\[
\lim_{n\rightarrow\infty}\frac{1}{\mu(\Lambda_{n})}\int_{\Lambda_{n}}\left|
\omega\left(  \prod_{j=0}^{k}\tau_{g^{m_{j}}}(a)\right)  \right|
dg=\omega(a)^{k+1}>0.
\]
If $\left(  A,\omega,\tau,G\right)  $ is not weakly mixing, it has a
nontrivial compact factor $\left(  N,\omega,\tau,G\right)  $.
Continuity of $g\mapsto\tau_{g}(a)$ implies that
$g\mapsto\omega\left(  \prod_{j=0}^{k} \tau_{g^{m_{j}}}(a)\right)  $
and $g\mapsto\Vert\tau_{g^{m_{j}}} (a)-a\Vert_{\omega}$ are also
continuous. The Szemer\'{e}di property for $\left(
N,\omega,\tau,G\right)  $ then follows directly from Theorem 5.13.
\end{proof}

It is easily seen that for $G=\mathbb{Z}^{q}$ and $G=\mathbb{R}^{q}$
such $(\Lambda_{n})$ and $(\Gamma_{n})$ exist: For
$G=\mathbb{Z}^{q}$ let $\Lambda_{n}=\{-n,\ldots,n\}^{q}$ and
$\Gamma_{n}=\{-|m|n,-|m|n+1,\ldots,|m|n\}$. For $\mathbb{R}^{q}$,
take $\Lambda_{n}=[-n,n]^{q}$ and $\Gamma_{n}=|m|\Lambda_{n}$.

\begin{acknowledgment} We thank Richard de Beer, Willem Fouch\'{e},
Johan Swart and Gusti van Zyl for useful discussions. We also thank
the referee for carefully reading the manuscript.
\end{acknowledgment}

\end{document}